\documentclass[10.5pt]{article}
\usepackage[leqno]{amsmath}
\usepackage{amsfonts}
\usepackage{graphicx}

\usepackage{amsmath}
\usepackage{amssymb}
\usepackage{latexsym}
\usepackage{amsmath, amsfonts,amssymb, amsthm, euscript,makeidx,color,mathrsfs}

\def\sqr#1#2{{\vcenter{\vbox{\hrule height.#2pt
              \hbox{\vrule width.#2pt height#1pt \kern#1pt \vrule width.#2pt}
              \hrule height.#2pt}}}}
\def\signed #1{{\unskip\nobreak\hfil\penalty50
              \hskip2em\hbox{}\nobreak\hfil#1
              \parfillskip=0pt \finalhyphendemerits=0 \par}}
\def\endpf{\signed {$\sqr69$}}

\def\3n{\negthinspace \negthinspace \negthinspace }
\def\2n{\negthinspace \negthinspace }
\def\1n{\negthinspace }

\def\dbE{\mathbb{E}}
\def\dbF{\mathbb{F}}

\def\dbH{\mathbb{H}}

\def\dbP{\mathbb{P}}

\def\dbR{\mathbb{R}}
\def\dbS{\mathbb{S}}

\def\={\buildrel \triangle \over =}

\def\ds{\displaystyle}

\def\ns{\noalign{\ss}}
\def\a{\alpha}
\def\b{\beta}
\def\g{\gamma}
\def\d{\delta}
\def\e{\varepsilon}
\def\z{\zeta}

\def\l{\lambda}

\def\n{\nu}
\def\si{\sigma}

\def\f{\varphi}
\def\th{\theta}

\def\i{\infty}
\def\G{\Gamma}
\def\D{\Delta}
\def\Th{\Theta}
\def\L{\Lambda}

\def\F{\Phi}
\def\O{\Omega}

\def\cD{{\cal D}}

\def\cF{{\cal F}}

\def\cH{{\cal H}}

\def\cP{{\cal P}}

\def\cR{{\cal R}}

\def\cU{{\cal U}}

\def\cX{{\cal X}}

\def\cl{{\cal l}}
\def\no{\noindent}

\def\ss{\smallskip}
\def\ms{\medskip}

\def\q{\quad}
\def\qq{\qquad}
\def\hb{\hbox}

\def\liminf{\mathop{\underline{\rm lim}}}

\def\Ra{\mathop{\Rightarrow}}

\def\llan{\left\langle}
\def\rran{\right\rangle}
\def\lan{\big\langle}
\def\ran{\big\rangle}

\def\esssup{\mathop{\rm esssup}}

\def\h{\widehat}
\def\wt{\widetilde}

\def\cd{\cdot}

\def\ae{\hbox{\rm a.e.{ }}}
\def\as{\hbox{\rm a.s.{ }}}

\def\tr{\hbox{\rm tr$\,$}}

\def\cl{\overline}

\def\deq{\triangleq}
\def\les{\leqslant}
\def\ges{\geqslant}

\def\({\Big (}
\def\){\Big )}
\def\[{\Big[}
\def\]{\Big]}
\def\bde{\begin{definition}}
\def\ede{\end{definition}}
\def\be{\begin{equation}}
\def\bel{\begin{equation}\label}
\def\ee{\end{equation}}
\def\bt{\begin{theorem}}
\def\et{\end{theorem}}
\def\bc{\begin{corollary}}
\def\ec{\end{corollary}}
\def\bl{\begin{lemma}}
\def\el{\end{lemma}}
\def\bp{\begin{proposition}}
\def\ep{\end{proposition}}
\def\bas{\begin{assumption}}
\def\eas{\end{assumption}}
\def\br{\begin{remark}}
\def\er{\end{remark}}
\def\ba{\begin{array}}
\def\ea{\end{array}}
\def\ed{\end{document}}

\def\square#1{\vbox{\hrule\hbox{\vrule height#1%
     \kern#1\vrule}\hrule}}
\def\rectangle#1#2{\vbox{\hrule\hbox{\vrule height#1%
     \kern#2\vrule}\hrule}}

\font\tenbb=msbm10 \font\sevenbb=msbm7 \font\fivebb=msbm5

\newfam\bbfam
\scriptscriptfont\bbfam=\fivebb \textfont\bbfam=\tenbb
\scriptfont\bbfam=\sevenbb

\newtheorem{lemma}{\hskip 1.4em Lemma}[section]
\newtheorem{remark}{\hskip 1.4em Remark}[section]

\newtheorem{theorem}{\hskip 1.4em Theorem}[section]
\newtheorem{corollary}{\hskip 1.4em Corollary}[section]

\newtheorem{definition}{\hskip 1.4em Definition}[section]
\newtheorem{proposition}{\hskip 1.4em Proposition}[section]
\newtheorem{assumption}{\hskip 1.4em Assumption}[section]

\makeatletter
   
   \@addtoreset{equation}{section}
\makeatother

\begin{document}

\title{\bf Open-Loop and Closed-Loop Solvabilities for Stochastic Linear Quadratic Optimal Control
Problems\thanks{This work is supported in part by NSF Grant
DMS-1406776.}}
\author{Jingrui Sun\thanks{Department of Applied Mathematics, The Hong Kong Polytechnic University, Hong Kong, China
(sjr@mail.ustc.edu.cn).}\ , \ Xun Li\thanks{Department of
Applied Mathematics, The Hong Kong Polytechnic University, Hong
Kong, China (malixun@polyu.edu.hk).}\ , \ and \ Jiongmin
Yong\thanks{Department of Mathematics, University of Central
Florida, Orlando, FL 32816, USA (jiongmin.yong@ucf.edu).}
 }
\maketitle

\noindent {\bf Abstract:} This paper is concerned with a stochastic
linear quadratic (LQ, for short) optimal control problem. The
notions of open-loop and closed-loop solvabilities are introduced. A
simple example shows that these two solvabilities are different.
Closed-loop solvability is established by means of solvability of
the corresponding Riccati equation, which is implied by the uniform
convexity of the quadratic cost functional. Conditions ensuring the
convexity of the cost functional are discussed, including the issue
that how negative the control weighting matrix-valued function
$R(\cd)$ can be. Finiteness of the LQ problem is characterized by
the convergence of the solutions to a family of Riccati equations.
Then, a minimizing sequence, whose convergence is equivalent to the
open-loop solvability of the problem, is constructed. Finally, an
illustrative example is presented.

\ms


\ms

\no\bf AMS Mathematics Subject Classification. \rm 49N10, 49N35, 93E20.

\section{Introduction}

Let $(\O,\cF,\dbF,\dbP)$ be a complete filtered probability space on
which a standard one-dimensional Brownian motion $W=\{W(t); 0\les t
< \i \}$ is defined, where $\dbF=\{\cF_t\}_{t\ges0}$ is the natural
filtration of $W$ augmented by all the $\dbP$-null sets in $\cF$
\cite{Karatzas-Shreve 1991,Yong-Zhou 1999}. Consider the following
controlled linear stochastic differential equation (SDE, for short)
on a finite horizon $[t,T]$:
\bel{state}\left\{\2n\ba{ll}
\ns\ds dX(s)=\big[A(s)X(s)+B(s)u(s)+b(s)\big]ds+\big[C(s)X(s)+D(s)u(s)+\si(s)\big]dW(s),
\qq s\in[t,T], \\
\ns\ds X(t)=x,\ea\right.\ee
where $A(\cd), B(\cd), C(\cd), D(\cd)$ are given deterministic
matrix-valued functions of proper dimensions, and $b(\cd), \si(\cd)$
are vector-valued $\dbF$-progressively measurable processes. In the
above, $X(\cd)$, valued in $\dbR^n$, is the {\it state process}, and
$u(\cd)$, valued in $\dbR^m$, is the {\it control process}. For any
$t\in[0,T)$, we introduce the following Hilbert space:
$$\cU[t,T]=\left\{u:[t,T]\times\O\to\dbR^m\bigm|u(\cd)\hb{ is $\dbF$-progressively
measurable, }\dbE\int_t^T|u(s)|^2ds<\i\right\}.$$
Any $u(\cd)\in\cU[t,T]$ is called an {\it admissible control} (on
$[t,T]$). Under some mild conditions on the coefficients, for any
{\it initial pair} $(t,x)\in[0,T)\times\dbR^n$ and admissible
control $u(\cd)\in\cU[t,T]$, (\ref{state}) admits a unique strong
solution $X(\cd)\equiv X(\cd\,;t,x,u(\cd))$. Next we introduce the
following cost functional:
\bel{cost}\ba{ll}
\ns\ds J(t,x;u(\cd))\deq\dbE\Bigg\{\lan GX(T),X(T)\ran+2\lan g,X(T)\ran\\
\ns\ds\qq\qq\qq\qq+\int_t^T\left[\llan\begin{pmatrix}Q(s)&\1nS(s)^\top\\S(s)&\1nR(s)\end{pmatrix}
                                    \begin{pmatrix}X(s)\\ u(s)\end{pmatrix},
                                    \begin{pmatrix}X(s)\\u(s)\end{pmatrix}\rran
+2\llan\begin{pmatrix}q(s)\\ \rho(s)\end{pmatrix},\begin{pmatrix}X(s)\\u(s)\end{pmatrix}\rran\right] ds\Bigg\},\ea\ee
where $G$ is a symmetric matrix, $Q(\cd)$, $S(\cd)$, $R(\cd)$ are deterministic matrix-valued
functions of proper dimensions with $Q(\cd)^\top=Q(\cd)$, $R(\cd)^\top=R(\cd)$; $g$ is allowed
to be an $\cF_T$-measurable random variable and $q(\cd), \rho(\cd)$
are allowed to be vector-valued $\dbF$-progressively measurable processes.
The classical stochastic LQ optimal control problem can be stated as follows.

\ms

\bf Problem (SLQ). \rm For any given initial pair
$(t,x)\in[0,T)\times\dbR^n$, find a $u^*(\cd)\in\cU[t,T]$, such that
\bel{optim}J(t,x;u^*(\cd))=\inf_{u(\cd)\in\cU[t,T]}J(t,x;u(\cd))\deq V(t,x).\ee

It is well-accepted that any $u^*(\cd)\in\cU[t,T]$ satisfying
(\ref{optim}) is called an {\it optimal control} of Problem (SLQ)
for the initial pair $(t,x)$, and the corresponding $X^*(\cd)\equiv
X(\cd\,; t,x,u^*(\cd))$ is called an {\it optimal state process};
the pair $(X^*(\cd),u^*(\cd))$ is called an {\it optimal pair}. The
function $V(\cd\,,\cd)$ is called the {\it value function} of
Problem (SLQ). When $b(\cd), \si(\cd), g, q(\cd), \rho(\cd)=0$, we denote the corresponding Problem (SLQ) by Problem $\hb{(SLQ)}^0$.
The corresponding cost functional and value function are denoted by $J^0(t,x;u(\cd))$ and $V^0(t,x)$, respectively.

\ms

When the stochastic part is absent, with $b(\cd)$, $g$, $q(\cd)$ and
$\rho(\cd)$ being deterministic, Problem (SLQ) is reduced to a
deterministic LQ optimal control problem, called Problem (DLQ).
Hence, Problem (DLQ) can be regarded as a special case of Problem
(SLQ). The history of Problem (DLQ) can be traced back to the work
of Bellman--Glicksberg--Gross \cite{Belman-Gicksberg-Gross 1958} in
1958, Kalman \cite{Kalman 1960} in 1960 and Letov \cite{Letov 1961}
in 1961. In the deterministic case, it is well known that
$R(s)\ges0$ is necessary for Problem (DLQ) to be finite (meaning
that the infimum of the cost functional is finite). When the control
weighting matrix $R(s)$ in the cost function is uniformly positive
definite, under some mild additional conditions on the other
weighting coefficients, the problem can be solved elegantly via the
Riccati equation; see \cite{Anderson-Moore 1989} for a thorough
study of the Riccati equation approach (see also \cite{Yong-Zhou
1999}). Stochastic LQ problems was firstly studied by Wonham
\cite{Wonham 1968} in 1968, followed by several researchers later
(see \cite{Davis 1977, Bensoussan 1983}, for examples). In those
works, the assumption that $R(s)>0$ was taken for granted. More
precisely, under the following {\it standard} conditions:
\bel{classical}G\ges0,\qq R(s)\ges\d I,\qq Q(s)-S(s)^\top
R(s)^{-1}S(s)\ges0,\qq\ae~s\in[0,T],\ee
for some $\d>0$, the corresponding Riccati equation is uniquely
solvable and Problem (SLQ) admits a unique optimal control which has
a linear state feedback representation (see \cite[Chapter
6]{Yong-Zhou 1999} or \cite{Chen-Zhou 2000}). In 1998,
Chen--Li--Zhou \cite{Chen-Li-Zhou 1998} found that Problem (SLQ)
might still be solvable even if $R(s)$ is not positive
semi-definite. See also some follow-up works of Lim--Zhou
\cite{Lim-Zhou 1999}, Chen--Zhou \cite{Chen-Zhou 2000}, and
Chen--Yong \cite{Chen-Yong 2001}, as well as the works of Hu--Zhou
\cite{Hu-Zhou 2003} and Qian--Zhou \cite{Qian-Zhou 2013} on the
study of solvability of indefinite Riccati equations (under certain
technical conditions). In 2001, Ait Rami--Moore--Zhou \cite{Ait
Rami-Moore-Zhou 2001} introduced a generalized Riccati equation
involving the pseudo-inverse of a matrix and an additional algebraic
constraint; see also Ait Rami--Zhou--Moore \cite{Rami-Zhou-Moore 2000}
for stochastic LQ optimal control problem on $[0,\i)$, and a
follow-up work of Wu--Zhou \cite{Wu-Zhou 2001}. Recently, based on
the work of Yong \cite{Yong 2013}, Huang--Li--Yong
\cite{Huang-Li-Yong 2014} studied a mean-field LQ optimal control
problem on $[0,\i)$. For stochastic LQ optimal control problems with
random coefficients, we further refer to the works of Chen--Yong
\cite{Chen-Yong 2000}, Tang \cite{Tang 2003}, and Kohlmann--Tang
\cite{Kohlmann-Tang 2003}.

\ms

Most recently, Sun--Yong \cite{Sun-Yong 2014} established that the
existence of open-loop optimal controls is equivalent to the
solvability of the corresponding optimality system which is a
forward-backward stochastic differential equation (FBSDE, for
short), and the existence of closed-loop optimal strategies is
equivalent to the existence of a regular solution to the
corresponding Riccati equation. From this point of view, this paper
can be regarded as a continuation of \cite{Sun-Yong 2014}, in a
certain sense. Inspired by a result found in \cite{Yong 2015}, we
are able to cook up an example for which open-loop optimal controls
exist, but the closed-loop optimal strategy does not exist. Because
of this, it is necessary to distinguish open-loop and closed-loop
solvabilities of Problem (SLQ). Next, having the equivalence between
the solvability of the Riccati equation and the closed-loop
solvability of Problem (SLQ), it is natural to seek conditions under
which the Riccati equation is solvable, and the sought conditions
are expected to be more general than (\ref{classical}) so that they
could include some (although might not be all) cases that $R(\cd)$
is allowed to be not positive semi-definite. One of our main results
in this paper is to establish the equivalence between the {\it
strongly regular solvability} of the Riccati equation (see below for
a definition) and the uniform convexity of the cost functional. Note
that uniform convexity condition is much weaker than
(\ref{classical}), and is different from conditions imposed in
\cite{Qian-Zhou 2013}.

\ms

Finiteness of Problem (SLQ) (meaning that the infimum of the cost
functional is finite) is another important issue. The notion was
introduced in \cite{Yong-Zhou 1999} (see also \cite{Chen-Yong
2000,Chen-Yong 2001}). Some investigations were carried out in
\cite{Mou-Yong 2006}. In this paper, due to the perfect structure of
Problem $\hb{(SLQ)}^0$, its finiteness will be characterized by the
convergence of the solutions to a family of Riccati equations.
As a by-product, we will construct minimizing sequences of Problem (SLQ)
in a very natural way, and the convergence of the sequences will lead to
the open-loop solvability of Problem (SLQ).

\ms

Among other things, we find several interesting facts which are
listed below:

\ms

\it Fact \rm 1. The value function $V(t,x)$ is not necessarily
continuous in $t$ even if Problem (SLQ) admits a {\it continuous}
open-loop optimal control at all initial pair $(t,x)\in[0,T)\times\dbR^n$.

\ms

\it Fact \rm 2. If Problem $\hb{(SLQ)}^0$ is finite at $t$, then it is finite at any
$t^\prime>t$.

\ms

\it Fact \rm 3. For Problem (SLQ) with $D(\cd)=0$, under the
assumption that $R(\cd)$ is uniformly positive definite, without
assuming the non-negativity of $Q(\cd)$ and $G$, the finiteness and
the unique closed-loop solvability of Problem (SLQ) are equivalent, which are
also equivalent to the uniform convexity of the cost functional.

\ms

\it Fact \rm 4. The existence of a regular solution to the Riccati
equation implies the open-loop solvability of Problem {\rm(SLQ)}.
However, it may happen that for any initial pair
$(t,x)\in[0,T)\times\dbR^n$, Problem {\rm(SLQ)} admits a {\it
continuous} open-loop optimal control, while the Riccati equation
does not admit a regular solution. This corrects an incorrect result
found in \cite{Ait Rami-Moore-Zhou 2001} (see Section 4 for
details).

\ms

The rest of the paper is organized as follows. Section 2 collects
some preliminary results. In Section 3, we study the cost functional
from a Hilbert space viewpoint and represent it as a quadratic
functional of $u(\cd)$. Section 4 is devoted to the strongly regular
solvability of the Riccati equation under the uniform convexity of
the cost functional. In Section 5, we discuss the finiteness of
Problem (SLQ) as well as the convexity of the cost functional. In
Section 6, we characterize the open-loop solvability of Problem
(SLQ) by means of the convergence of minimizing sequences. An
example is presented in section 7 to illustrate some relevant
results obtained.

\section{Preliminaries}

We recall that $\dbR^n$ is the $n$-dimensional Euclidean space,
$\dbR^{n\times m}$ is the space of all $(n\times m)$ matrices,
endowed with the inner product $\langle M,N\rangle\mapsto\tr[M^\top
N]$ and the norm $|M|=\sqrt{\tr[M^\top M]}$, $\dbS^n\subseteq
\dbR^{n\times n}$ is the set of all $(n\times n)$ symmetric
matrices, $\cl{\dbS^n_+}\subseteq\dbS^n$ is the set of all $(n\times
n)$  positive semi-definite matrices, and
$\dbS^n_+\subseteq\cl{\dbS^n_+}$ is the set of all $(n\times n)$
positive-definite matrices. When there is no confusion, we shall use
$\langle \cd\,,\cd\rangle$ for inner products in possibly different
Hilbert spaces. Also, $M^\dag$ stands for the (Moore-Penrose)
pseudo-inverse of the matrix $M$ (\cite{Penrose 1955}), and $\cR(M)$
stands for the range of the matrix $M$. Next, let $T>0$ be a fixed
time horizon. For any $t\in[0,T)$ and Euclidean space $\dbH$, let
$$\ba{ll}
C([t,T];\dbH)=\Big\{\f:[t,T]\to\dbH\bigm|\f(\cd)\hb{ is
continuous }\1n\Big\},\\
\ns\ds
L^p(t,T;\dbH)=\left\{\f:[t,T]\to\dbH\biggm|\int_t^T|\f(s)|^pds<\i\right\},\q1\les p<\i,\\
\ns\ds
L^\infty(t,T;\dbH)=\left\{\f:[t,T]\to\dbH\biggm|\esssup_{s\in[t,T]}|\f(s)|<\i\right\}.\ea$$
We denote
$$\ba{ll}
\ns\ds L^2_{\cF_T}(\O;\dbH)=\Big\{\xi:\O\to\dbH\bigm|\xi\hb{ is
$\cF_T$-measurable, }\dbE|\xi|^2<\i\Big\},\\
\ns\ds
L_\dbF^2(t,T;\dbH)=\left\{\f:[t,T]\times\O\to\dbH\bigm|\f(\cd)\hb{ is
$\dbF$-progressively measurable},\dbE\int^T_t|\f(s)|^2ds<\i\right\},\\
\ns\ds
L_\dbF^2(\O;C([t,T];\dbH))=\left\{\f:[t,T]\times\O\to\dbH\bigm|\f(\cd)\hb{
is $\dbF$-adapted, continuous, }\dbE\left[\sup_{s\in[t,T]}|\f(s)|^2\right]<\i\right\},\\
\ns\ds L^2_\dbF(\O;L^1(t,T;\dbH))=\left\{\f:[t,T]\times
\O\to\dbH\bigm|\f(\cd)\hb{ is $\dbF$-progressively measurable},
\dbE\left(\int_t^T|\f(s)|ds\right)^2<\i\right\}.\ea$$
For an $\dbS^n$-valued function $F(\cd)$ on $[t,T]$, we use the
notation $F(\cd)\gg0$ to indicate that $F(\cd)$ is uniformly
positive definite on $[t,T]$, i.e., there exists a constant
$\delta>0$ such that
$$F(s)\ges\delta I,\qq\ae~s\in[t,T].$$

The following standard assumptions will be in force throughout this paper.

\ms

{\bf(H1)} The coefficients of the state equation satisfy the following:
$$\left\{\2n\ba{ll}
\ns\ds A(\cd)\in L^1(0,T;\dbR^{n\times n}),\q B(\cd)\in L^2(0,T;\dbR^{n\times m}), \q b(\cd)\in L^2_\dbF(\O;L^1(0,T;\dbR^n)),\\
\ns\ds C(\cd)\in L^2(0,T;\dbR^{n\times n}),\q D(\cd)\in L^\i(0,T;\dbR^{n\times m}), \q\si(\cd)\in L_\dbF^2(0,T;\dbR^n).\ea\right.$$

{\bf(H2)} The weighting coefficients in the cost functional satisfy the following:
$$\left\{\2n\ba{ll}
\ns\ds  G\in\dbS^n,\q Q(\cd)\in L^1(0,T;\dbS^n),\q S(\cd)\in L^2(0,T;\dbR^{m\times
n}),\q R(\cd)\in L^\infty(0,T;\dbS^m),\\
\ns\ds g\in L^2_{\cF_T}(\O;\dbR^n),\q q(\cd)\in L^2_\dbF(\O;L^1(0,T;\dbR^n)),\q\rho(\cd)\in
L_\dbF^2(0,T;\dbR^m).\ea\right.$$

By \cite[Proposition 2.1]{Sun-Yong 2014}, under (H1)--(H2), for any
$(t,x)\in[0,T)\times\dbR^n$ and $u(\cd)\in\cU[t,T]$, the state
equation (\ref{state}) admits a unique strong solution $X(\cd)\equiv
X(\cd\,;t,x,u(\cd))$, and the cost functional (\ref{cost}) is
well-defined. Then Problem (SLQ) makes sense. It is worthy of
pointing out that in (H2), we do not impose any
positive-definiteness/non-negativeness conditions on the weighting
matrix/matrix-valued functions $G$, $Q(\cd)$ and $R(\cd)$. We now
introduce the following definition.

\ms

\bf Definition 2.1. \rm (i) Problem (SLQ) is said to be {\it finite
at initial pair $(t,x)\in[0,T]\times\dbR^n$} if
\bel{V>-infty}V(t,x)>-\i.\ee
Problem (SLQ) is said to be {\it finite at $t\in[0,T]$} if
(\ref{V>-infty}) holds for all $x\in\dbR^n$, and Problem (SLQ) is
said to be {\it finite} if (\ref{V>-infty}) holds for all
$(t,x)\in[0,T]\times\dbR^n$.

\ms

(ii) An element $u^*(\cd)\in\cU[t,T]$ is called an {\it open-loop
optimal control} of Problem (SLQ) for the initial pair
$(t,x)\in[0,T]\times\dbR^n$ if
\bel{open-opti}J(t,x;u^*(\cd))\les J(t,x;u(\cd)),\qq\forall
u(\cd)\in\cU[t,T].\ee
If an open-loop optimal control (uniquely) exists for
$(t,x)\in[0,T]\times\dbR^n$, Problem (SLQ) is said to be ({\it
uniquely}) {\it open-loop solvable at $(t,x)\in[0,T]\times\dbR^n$}.
Problem (SLQ) is said to be ({\it uniquely}) {\it open-loop solvable
at $t\in[0,T)$} if for the given $t$, (\ref{open-opti}) holds for
all $x\in\dbR^n$, and Problem (SLQ) is said to be ({\it uniquely})
{\it open-loop solvable} (on $[0,T)\times\dbR^n$) if it is
(uniquely) open-loop solvable at all $(t,x)\in[0,T)\times\dbR^n$.

\ms

(iii) A pair $(\Th^*(\cd),v^*(\cd))\in L^2(t,T;\dbR^{m\times
n})\times\cU[t,T]$ is called a {\it closed-loop optimal strategy} of
Problem (SLQ) on $[t,T]$ if
\bel{closed-opti}\ba{ll}
\ns\ds J(t,x;\Th^*(\cd)X^*(\cd)+v^*(\cd))\les
J(t,x;u(\cd)),\qq\forall x\in\dbR^n,\q u(\cd)\in\cU[t,T],\ea\ee
where $X^*(\cd)$ is the strong solution to the following closed-loop
system:
\bel{closed-loop0}\left\{\2n\ba{ll}
\ns\ds dX^*(s)=\Big\{\big[A(s)+B(s)\Th^*(s)\big]X^*(s)+B(s)v^*(s)+b(s)\Big\}ds\\
\ns\ds\qq\qq\qq+\Big\{\big[C(s)+D(s)\Th^*(s)\big]X^*(s)+D(s)v^*(s)+\si(s)\Big\}dW(s),
\qq s\in[t,T],\\
\ns\ds X^*(t)=x.\ea\right.\ee
If a closed-loop optimal strategy (uniquely) exists on $[t,T]$,
Problem (SLQ) is said to be ({\it uniquely}) {\it closed-loop
solvable on $[t,T]$}. Problem (SLQ) is said to be ({\it uniquely})
{\it closed-loop solvable} if it is (uniquely) closed-loop solvable
on any $[t,T]$.

\ms

\rm

We emphasize that in the definition of closed-loop optimal strategy,
(\ref{closed-opti}) has to be true for all $x\in\dbR^n$. One sees
that if $(\Th^*(\cd),v^*(\cd))$ is a closed-loop optimal strategy of
problem (SLQ) on $[t,T]$, then the outcome
$u^*(\cd)\equiv\Th^*(\cd)X^*(\cd)+v^*(\cd)$ is an open-loop optimal
control of Problem (SLQ) for the initial pair $(t,X^*(t))$. Hence,
the existence of closed-loop optimal strategies implies the
existence of open-loop optimal controls. But, the existence of
open-loop optimal controls does not necessarily imply the existence
of a closed-loop optimal strategy. Here is such an example which is
inspired by an example for deterministic LQ problems found in
\cite{Yong 2015}.

\ms

\bf Example 2.2. \rm Consider the following Problem $\hb{(SLQ)}^0$ with
one-dimensional state equation:
\bel{ex2.2-1}\left\{\2n\ba{ll}
\ns\ds dX(s)=\big[u_1(s)+u_2(s)\big]ds+\big[u_1(s)-u_2(s)\big]dW(s),\qq s\in[t,1], \\
\ns\ds X(t)=x,\ea\right.\ee
and cost functional
\bel{ex2.2-2}J^0(t,x;u(\cd))=\dbE X(1)^2.\ee
In this example, $u(\cd)=(u_1(\cd),u_2(\cd))^\top$. It is clear that
$$V^0(t,x)=\inf_{u(\cd)\in\cU[t,T]}J^0(t,x;u(\cd))\ges0,\qq\forall(t,x)\in[0,1]\times\dbR.$$
On the other hand, for any $t\in[0,1)$, $\b\ges{1\over 1-t}$, by
taking
$$u^\b(s)=-{\b x\over2}{\bf1}_{[t,t+{1\over\b}]}(s)\begin{pmatrix}1\\ 1\end{pmatrix},\qq
s\in[t,1],$$
we have
$$X(s)=0,\qq s\in[t+1/\b,1].$$
Hence,
$$J^0(t,x;u^\b(\cd))=0,\qq(t,x)\in[0,1)\times\dbR.$$
This implies that $\big\{u^\b(\cd)\bigm|\b\ges{1\over1-t}\big\}$ is
a family of open-loop optimal controls of the corresponding Problem
$\hb{(SLQ)}^0$ for the initial pair $(t,x)\in[0,1)\times\dbR$, and
therefore,
$$V^0(t,x)=\left\{\2n\ba{ll}
\ns\ds0,\qq\ (t,x)\in[0,1)\times\dbR,\\
\ns\ds x^2,\qq t=1,~x\in\dbR,\ea\right.$$
which is discontinuous at $t=1$, $x\ne0$. Note also that if we take
$\b={1\over1-t}$, then the corresponding open-loop optimal control,
denoted by $\bar u(\cd)$, is given by
\bel{bar u}\bar u(s)=-{x\over2(1-t)}\begin{pmatrix}1\\
1\end{pmatrix},\qq s\in[t,1],\ee
which is a constant vector (only depends on the initial pair
$(t,x)$). Thus, it is continuous (in $s\in[t,1]$). Now, we claim
that this Problem $\hb{(SLQ)}^0$ is not closed-loop solvable on any
$[t,1]$ with $t\in[0,1)$. In fact, if for some $t\in[0,1)$, there
exists a closed-loop optimal strategy $(\Th^*(\cd),v^*(\cd))$, then
by definition, one must have
\bel{ex2.2-3}0\les J^0(t,x;\Th^*(\cd)X^*(\cd)+v^*(\cd))\les J^0(t,x;u^\b(\cd))=0,\qq\forall x\in\dbR.\ee
Let
$$\Th^*(\cd)=\begin{pmatrix}\Th_1^*(\cd)\\
\Th^*_2(\cd)\end{pmatrix},\qq v^*(\cd)=\begin{pmatrix}v^*_1(\cd)\\
v^*_2(\cd)\end{pmatrix}.$$
Then we have from (\ref{ex2.2-3}) that the solution $X^*(\cd)$ of the following
closed-loop system:
$$\left\{\2n\ba{ll}
\ns\ds dX^*(s)=\Big\{\big[\Th^*_1(s)+\Th_2^*(s)\big]X^*(s)+\big[v^*_1(s)+v^*_2(s)\big]\Big\}ds\\
\ns\ds\qq\qq~~+\Big\{\big[\Th^*_1(s)-\Th_2^*(s)\big]X^*(s)+\big[v^*_1(s)-v^*_2(s)\big]\Big\}dW(s),
\qq s\in[t,1],\\
\ns\ds X^*(t)=x,\ea\right.$$
must satisfy
$$X^*(1)=0,\qq\forall x\in\dbR.$$
Note that
$$\left\{\2n\ba{ll}
\ns\ds d\big[\dbE X^*(s)\big]=\Big\{\big[\Th^*_1(s)+\Th_2^*(s)\big]\dbE X^*(s)+\dbE\big[v^*_1(s)+v^*_2(s)\big]\Big\}ds,\qq
s\in[t,1],\\
\ns\ds \dbE X^*(t)=x.\ea\right.$$
Consequently,
$$\ba{ll}
\ns\ds0=\dbE
X^*(1)=e^{\int_t^1[\Th_1^*(s)+\Th^*_2(s)]ds}x+\int_t^1e^{\int_r^1[\Th^*_1(s)
+\Th^*_2(s)]ds}\dbE[v_1^*(r)+v^*_2(r)]dr,\qq\forall x\in\dbR,\ea$$
which is impossible since the above has to be true for all
$x\in\dbR$. Hence, the corresponding Problem $\hb{(SLQ)}^0$ is not
closed-loop solvable on any $[t,1]$ with $t\in[0,1)$, although the
problem admits a continuous open-loop optimal control for any
initial pair $(t,x)\in[0,1)\times\dbR^n$.

\ms

Due to the above indicated situation, unlike in \cite{Sun-Yong
2014}, and in classical literature on LQ problems, we distinguish
the notions of open-loop and closed-loop solvabilities for Problem
(SLQ). We repeat here that for given initial time $t\in[0,T)$, an
open-loop optimal control is allowed to depend on the initial state
$x$, whereas, a closed-loop optimal strategy is required to be
independent of the initial state $x$. Because of the nature of
closed-loop strategies, we define the finiteness of Problem (SLQ)
only in terms of open-loop controls.

\ms

To conclude this section, we present some lemmas which will be used
frequently in sequel.

\ms

\bf Lemma 2.3. \sl Let {\rm(H1)--(H2)} hold and $\Th(\cd)\in
L^2(0,T;\dbR^{m\times n})$. Let $P(\cd)\in C([0,T];\dbS^n)$ be the
solution to the following Lyapunov equation:
\bel{P-Th}\left\{\2n\ba{ll}
\ns\ds\dot P+P(A+B\Th)+(A+B\Th)^\top P+(C+D\Th)^\top P(C+D\Th)\\
\ns\ds\q+\,\Th^\top R\Th+S^\top\Th+\Th^\top S+Q=0,\qq\ae~s\in[0,T],\\
\ns\ds P(T)=G.\ea\right.\ee
Then for any $(t,x)\in[0,T)\times\dbR^n$ and $u(\cd)\in\cU[t,T]$, we have
$$\ba{ll}
\ns\ds J^0(t,x;\Th(\cd)X(\cd)+u(\cd))=\langle P(t)x,x\rangle+\dbE\int_t^T\Big\{2\lan\big[B^\top P+D^\top PC+S+(R+D^\top PD)\Th\big]X,u\ran\\
\ns\ds\qq\qq\qq\qq\qq\qq\qq\qq\qq\qq~\ +\lan (R+D^\top PD)u,u\ran\Big\}ds.
\ea$$

\ms

\it Proof. \rm For any $(t,x)\in[0,T)\times\dbR^n$ and $u(\cd)\in\cU[t,T]$,
let $X(\cd)$ be the solution of
$$\left\{\2n\ba{ll}
\ns\ds dX(s)=\big[(A+B\Th)X+Bu\big]ds+\big[(C+D\Th)X+Du\big]dW(s),\qq s\in[t,T], \\
\ns\ds X(t)=x.\ea\right.$$
Applying It\^o's formula to $s\mapsto\langle P(s)X(s),X(s)\rangle$, we have
$$\ba{ll}
\ns\ds J^0(t,x;\Th(\cd)X(\cd)+u(\cd))
=\dbE\left\{\langle GX(T),X(T)\rangle
+\int_t^T\llan\begin{pmatrix}Q&S^\top\\S&R\end{pmatrix}
             \begin{pmatrix}X\\ \Th X+u\end{pmatrix},
             \begin{pmatrix}X\\ \Th X+u\end{pmatrix}\rran ds\right\}\\
\ns\ds=\langle P(t)x,x\rangle+\dbE\int_t^T\Big\{\lan\dot PX,X\ran+\lan P\big[(A+B\Th)X+Bu\big],X\ran+\lan PX,(A+B\Th)X+Bu\ran\\
\ns\ds\qq\qq\qq\qq\qq~+\lan P\big[(C+D\Th)X+Du\big],(C+D\Th)X+Du\ran+\lan QX,X\ran\\
\ns\ds\qq\qq\qq\qq\qq~+2\lan SX,\Th X+u\ran+\lan R(\Th X+u),\Th X+u\ran\Big\}ds\\
\ns\ds=\langle P(t)x,x\rangle+\dbE\int_t^T\Big\{\lan\big[\dot P+P(A+B\Th)+(A+B\Th)^\top P+(C+D\Th)^\top P(C+D\Th)\\
\ns\ds\qq\qq\qq\qq\qq\q~~\1n+Q+\Th^\top S+S^\top\Th+\Th^\top R\Th\big]X,X\ran\\
\ns\ds\qq\qq\qq\qq\qq~+2\lan\big[B^\top P+D^\top PC+S+(R+D^\top PD)\Th\big]X,u\ran+\lan (R+D^\top PD)u,u\ran\Big\}ds\\
\ns\ds=\langle P(t)x,x\rangle+\dbE\int_t^T\Big\{2\lan\big[B^\top P+D^\top PC+S+(R+D^\top PD)\Th\big]X,u\ran+\lan (R+D^\top PD)u,u\ran\Big\}ds.
\ea$$
This completes the proof. \endpf

\ms

The following lemma is concerned with the solution to a Lyapunov
equation, whose proof can be found in \cite{Chen-Zhou 2000} (see
also \cite[Chapter 6]{Yong-Zhou 1999}).

\ms

\bf Lemma 2.4. \sl Let
$$\wt A(\cd)\in L^1(0,T;\dbR^{n \times n}),\q \wt C(\cd)\in L^2(0,T;\dbR^{n \times n}),
\q\wt Q(\cd)\in L^1(0,T;\dbS^n), \q \wt G\in\dbS^n.$$
Then the following Lyapunov equation
$$\left\{\2n\ba{ll}
\ns\ds \dot{P}(s)+P(s)\wt A(s)+\wt A(s)^\top P(s)+\wt C(s)^\top P(s)\wt C(s)+\wt Q(s)=0,\qq\ae~s\in [0,T],\\
\ns\ds P(T)=\wt G,\ea\right.$$
admits a unique solution $P(\cd) \in C([0,T];\dbS^n)$ given by
$$P(t)=\dbE\left\{\big[\wt \F(T)\wt \F(t)^{-1}\big]^\top \wt G\big[\wt \F(T)\wt \F(t)^{-1}\big]
+\int_t^T\big[\wt \F(s)\wt \F(t)^{-1}\big]^\top \wt Q(s)\big[\wt \F(s)\wt \F(t)^{-1}\big]ds\right\},$$
where $\wt \F(\cd)$ is the solution of
$$\left\{\2n\ba{ll}
\ns\ds d\wt \F(s)=\wt A(s)\wt \F(s)ds+\wt C(s)\wt \F(s)dW(s),\qq s\ges0,\\
\ns\ds\wt \F(0)=I.\ea\right.$$
Consequently, if
$$\wt G\ges0,\qq\wt Q(s)\ges0,\qq \ae~s\in[0,T],$$
then $P(\cd)\in C([0,T];\cl{\dbS^n_+})$.

\ms

\bf Lemma 2.5. \sl For any $u(\cd)\in\cU[t,T]$, let $X^{(u)}(\cd)$
be the solution of
\bel{}\left\{\2n\ba{ll}
\ns\ds dX^{(u)}(s)=\big[A(s)X^{(u)}(s)+B(s)u(s)\big]ds+\big[C(s)X^{(u)}(s)+D(s)u(s)\big]dW(s),\qq s\in[t,T], \\
\ns\ds X^{(u)}(t)=0.\ea\right.\ee
Then for any $\Th(\cd)\in L^2(t,T;\dbR^{m\times n})$, there exists a constant $\gamma>0$ such that
\bel{lem-2.6}\dbE\int_t^T\big|u(s)-\Th(s) X^{(u)}(s)\big|^2ds
\ges\g\dbE\int_t^T|u(s)|^2ds,\qq\forall u(\cd)\in\cU[t,T].\ee

\ms

\it Proof. \rm Let $\Th(\cd)\in L^2(t,T;\dbR^{m\times n})$. Define a bounded linear operator $\mathfrak{L}:\cU[t,T]\to\cU[t,T]$ by
$$\mathfrak{L}u=u-\Th X^{(u)}.$$
Then $\mathfrak{L}$ is bijective and its inverse $\mathfrak{L}^{-1}$ is given by
$$\mathfrak{L}^{-1}u=u+\Th\wt X^{(u)},$$
where $\wt X^{(u)}(\cd)$ is the solution of
$$\left\{\2n\ba{ll}
\ns\ds d\wt X^{(u)}(s)=\Big\{\big[A(s)+B(s)\Th(s)\big]\wt X^{(u)}(s)+B(s)u(s)\Big\}ds\\
\ns\ds\qq\qq\qq+\,\Big\{\big[C(s)+D(s)\Th(s)\big]\wt X^{(u)}(s)+D(s)u(s)\Big\}dW(s),\qq s\in[t,T], \\
\ns\ds\wt X^{(u)}(t)=0.\ea\right.$$
By the bounded inverse theorem, $\mathfrak{L}^{-1}$ is bounded with
norm $\|\mathfrak{L}^{-1}\|>0$. Thus,
$$\ba{ll}
\ns\ds\dbE\int_t^T|u(s)|^2ds=\dbE\int_t^T|(\mathfrak{L}^{-1}\mathfrak{L}u)(s)|^2ds
\les\|\mathfrak{L}^{-1}\|\dbE\int_t^T|(\mathfrak{L}u)(s)|^2ds\\
\ns\ds\qq\qq\qq~~\1n=\|\mathfrak{L}^{-1}\|\dbE\int_t^T\big|u(s)-\Th(s) X^{(u)}(s)\big|^2ds,
\qq\forall u(\cd)\in\cU[t,T],\ea$$
which implies (\ref{lem-2.6}) with $\g=\|\mathfrak{L}^{-1}\|^{-1}$.
\endpf

\ms

Finally, we state the following extended Schur's lemma whose proof
can be found in \cite{Albert 1969}.

\ms

\bf Lemma 2.6. \sl Let $Q\in\dbS^n$, $R\in\dbS^m$ and
$S\in\dbR^{m\times n}$. Then
\bel{2.11}\begin{pmatrix}Q&S^\top\\ S&R\end{pmatrix}\ges0,\ee
if and only if
\bel{2.12}Q-S^\top R^\dag S\ges0,\q R\ges0,\q\cR(S)\subseteq\cR(R).\ee

\rm

\ms

Note that the third condition in (\ref{2.12}) is equivalent to the
following:
\bel{2.13}S^\top(I-RR^\dag)=0.\ee

\section{Representation of the Cost Functional}

In this section, we will present a representation of the cost
functional for Problem (SLQ), from which we will obtain some basic
conditions ensuring the convexity of the cost functional. Convexity
of the cost functional will play a crucial role in the study of
finiteness and open-loop/closed-loop solvability of Problem (SLQ).
The following proposition is a summary of some relevant results
found in \cite{Yong-Zhou 1999}.

\ms

\bf Proposition 3.1. \sl Let {\rm(H1)--(H2)} hold. For any
$t\in[0,T)$, there exists a bounded self-adjoint linear operator
$M_2(t):\cU[t,T]\to\cU[t,T]$, a bounded linear operator
$M_1(t):\dbR^n\to\cU[t,T]$, an $M_0(t)\in\dbS^n$ and
$\nu_t(\cd)\in\cU[t,T], y_t\in\dbR^n, c_t\in\dbR$ such that
\bel{J-rep1}\ba{ll}
\ns\ds J(t,x;u(\cd))=\langle M_2(t)u,u\rangle+2\langle M_1(t)x,u\rangle+\langle M_0(t)x,x\rangle
+2\langle u,\nu_t\rangle+2\langle x,y_t\rangle+c_t,\\
\ns\ds J^0(t,x;u(\cd))=\langle M_2(t)u,u\rangle+2\langle M_1(t)x,u\rangle+\langle M_0(t)x,x\rangle,\\
\ns\ds\qq\qq\qq\qq\qq\qq\qq\qq\qq\forall(x,u(\cd))\in\dbR^n\times\cU[t,T].\ea\ee
Moreover, let $(X_0(\cd),Y(\cd),Z(\cd))$ be the adapted solution of
the following (decoupled) linear FBSDE:
\bel{FBSDE01}\left\{\2n\ba{ll}
\ns\ds
dX_0(s)=\big[A(s)X_0(s)+B(s)u(s)\big]ds+\big[C(s)X_0(s)+D(s)u(s)\big]dW(s),\qq
s\in[t,T],\\
\ns\ds dY(s)=-\big[A(s)^\top Y(s)+C(s)^\top
Z(s)+Q(s)X_0(s)+S(s)^\top u(s)\big]ds+Z(s)dW(s),\qq s\in[t,T],\\
\ns\ds X_0(t)=0,\qq Y(T)=GX_0(T).\ea\right.\ee
Then
\bel{L_2-rep}(M_2(t)u(\cd))(s)=B(s)^\top Y(s)+D(s)^\top
Z(s)+S(s)X_0(s)+R(s)u(s),\qq s\in[t,T].\ee
Let $(\bar X_0(\cd),\bar Y(\cd),\bar Z(\cd))$ be the adapted
solution to the following (decoupled) FBSDE:
\bel{FBSDE02}\left\{\2n\ba{ll}
\ns\ds d\bar X_0(s)=A(s)\bar X_0(s)ds+C(s)\bar X_0(s)dW(s),\qq
s\in[t,T],\\
\ns\ds d\bar Y(s)=-\big[A(s)^\top\bar Y(s)+C(s)^\top
\bar Z(s)+Q(s)\bar X_0(s)\big]ds+\bar Z(s)dW(s),\qq s\in[t,T],\\
\ns\ds\bar X_0(t)=x,\qq\bar Y(T)=G\bar X_0(T).\ea\right.\ee
Then
\bel{L_1-rep}\left\{\2n\ba{ll}
\ns\ds(M_1(t)x)(s)=B(s)^\top\bar Y(s)+D(s)^\top \bar Z(s)+S(s)\bar
X_0(s),\qq s\in[t,T],\\
\ns\ds M_0(t)x=\dbE\big[\bar Y(t)\big].\ea\right.\ee
Also, let $(\h X_0(\cd),\h Y(\cd),\h Z(\cd))$ be the adapted
solution to the following (decoupled) FBSDE:
\bel{FBSDE-nu}\left\{\2n\ba{ll}
\ns\ds d\h X_0(s)=\big[A(s)\h X_0(s)+b(s)\big]ds+\big[C(s)\h X_0(s)+\si(s)\big]dW(s),\qq
s\in[t,T],\\
\ns\ds d\h Y(s)=-\big[A(s)^\top\h Y(s)+C(s)^\top
\h Z(s)+Q(s)\h X_0(s)+q(s)\big]ds+\h Z(s)dW(s),\qq s\in[t,T],\\
\ns\ds\h X_0(t)=0,\qq\h Y(T)=G\h X_0(T)+g.\ea\right.\ee
Then
\bel{nu-rep}\nu_t(s)=B(s)^\top\h Y(s)+D(s)^\top \h Z(s)+S(s)\h X_0(s)+\rho(s),\qq s\in[t,T].\ee
Finally, $M_0(\cd)$ solves the following Lyapunov equation:
\bel{3.8}\left\{\2n\ba{ll}
\ns\ds\dot M_0(t)+M_0(t)A(t)+A(t)^\top M_0(t)+C(t)^\top
M_0(t)C(t)+Q(t)=0,\qq t\in[0,T],\\
\ns\ds M_0(T)=G,\ea\right.\ee
and it admits the following representation:
\bel{L_0}M_0(t)=\dbE\left\{\big[\F(T)\F(t)^{-1}\big]^\top
G\big[\F(T)\F(t)^{-1}\big]+\int_t^T\big[\F(s)\F(t)^{-1}\big]^\top
Q(s)\big[\F(s)\F(t)^{-1}\big]ds\right\},\ee
where $\F(\cd)$ is the solution to the following SDE for
$\dbR^{n\times n}$-valued process:
\bel{F}\left\{\2n\ba{ll}
\ns\ds d\F(s)=A(s)\F(s)ds+C(s)\F(s)dW(s),\qq s\ges0,\\
\ns\ds\F(0)=I.\ea\right.\ee

\ms

\it Proof. \rm Let $\F(\cd)$ be the solution to (\ref{F}). Then
$\F(s)^{-1}$ exists for all $s\ges0$ and the following holds:
$$\left\{\2n\ba{ll}
\ns\ds d\big[\F(s)^{-1}\big]=-\F(s)^{-1}\big[A(s)-C(s)^2\big]ds-\F(s)^{-1}C(s)dW(s),\qq s\ges0,\\
\ns\ds\F(0)^{-1}=I.\ea\right.$$
By the variation of constants formula, the solution $X(\cd)\equiv X(\cd\,;t,x,u(\cd))$ of the state
equation (\ref{state}) can be written as follows:
$$\ba{ll}
\ns\ds X(s)=\F(s)\F(t)^{-1}x+\F(s)\left\{\int_t^s\F(r)^{-1}\big[B(r)-C(r)D(r)\big]u(r)dr
+\int_t^s\F(r)^{-1}D(r)u(r)dW(r)\right\}\\
\ns\ds\qq\q~~+\F(s)\left\{\int_t^s\F(r)^{-1}\big[b(r)-C(r)\si(r)\big]dr+\int_t^s\F(r)^{-1}\si(r)dW(r)\right\},
\qq s\in[t,T].\ea$$
Now, let
$$h_t(s)=\F(s)\left\{\int_t^s\F(r)^{-1}\big[b(r)-C(r)\si(r)\big]dr+\int_t^s\F(r)^{-1}\si(r)dW(r)\right\},\qq s\in[t,T],$$
and define the following operators: For any $t\in[0,T)$,
$(x,u(\cd))\in\dbR^n\times\cU[t,T]$,
$$\left\{\2n\ba{ll}
\ns\ds(L_tu)(\cd)=\F(\cd)\left\{\int_t^\cd\F(r)^{-1}\big[B(r)-C(r)D(r)\big]u(r)dr+\int_t^\cd\F(r)^{-1}D(r)u(r)dW(r)\right\},\\
\ns\ds(\G_tx)(\cd)=\F(\cd)\F(t)^{-1}x,\q\h L_tu=(L_tu)(T),\q\h\G_tx=(\G_tx)(T).\ea\right.$$
Clearly, for any $t\in[0,T)$,
$$L_t:\cU[t,T]\to\cX[t,T],\qq\G_t:\dbR^n\to\cX[t,T],\qq\h L_t:\cU[t,T]\to\cX_T,
\qq\h\G_t:\dbR^n\to\cX_T$$
are all bounded linear operators, where $\cX[t,T]\equiv L_\dbF^2(t,T;\dbR^n)$
and $\cX_T\equiv L^2_{\cF_T}(\O;\dbR^n)$. Then, for
any $t\in[0,T)$ and $(x,u(\cd))\in\dbR^n\times\cU[t,T]$, the
corresponding state process $X(\cd)$ and its terminal value $X(T)$
are given by
$$\left\{\2n\ba{ll}
\ns\ds X(\cd)=(\G_tx)(\cd)+(L_tu)(\cd)+h_t(\cd),\\
\ns\ds X(T)=\h\G_tx+\h L_tu+h_t(T).\ea\right.$$
Hence, the cost functional can be written as
$$\ba{ll}
\ns\ds J(t,x;u(\cd))=\lan G\big(\h\G_tx+\h L_tu+h_t(T)\big),\h\G_tx+\h L_tu+h_t(T)\ran
+2\lan g,\h\G_tx+\h L_tu+h_t(T)\ran\\
\ns\ds\qq\qq\qq~+\lan Q(\G_tx+L_tu+h_t),\G_tx+L_tu+h_t\ran+2\lan S(\G_tx+L_tu+h_t),u\ran\\
\ns\ds\qq\qq\qq~+\langle Ru,u\rangle+2\langle q,\G_tx+L_tu+h_t\rangle+2\langle \rho,u\rangle.\ea$$
In the above, $\langle\cd\,,\cd\rangle$ is used for inner products in
possibly different spaces. Further, the adjoint operators
$$L_t^*:\cX[t,T]\to\cU[t,T],\qq\G_t^*:\cX[t,T]\to\dbR^n,\qq\h L_t^*:\cX_T\to\cU[t,T],\qq
\h\G_t^*:\cX_T\to\dbR^n$$
are given by the following:
$$(L^*_t\xi)(s)=B(s)^\top Y_0(s)+D(s)^\top Z_0(s),\q
s\in[t,T],\qq\q\G_t^*\xi=\dbE\big[Y_0(t)\big],$$
with $(Y_0(\cd),Z_0(\cd))$ being the adapted solution to the
following backward stochastic differential equation (BSDE, for
short):
\bel{BSDE0}\left\{\2n\ba{ll}
\ns\ds dY_0(s)=-\big[A(s)^\top Y_0(s)+C(s)^\top
Z_0(s)+\xi(s)\big]ds+Z_0(s)dW(s),\qq s\in[t,T],\\
\ns\ds Y_0(T)=0,\ea\right.\ee
and
$$(\h L_t^*\eta)(s)=B(s)^\top Y_1(s)+D(s)^\top Z_1(s),\q
s\in[t,T],\qq\h\G^*_t\eta=\dbE\big[Y_1(t)\big],$$
with $(Y_1(\cd),Z_1(\cd))$ being the adapted solution to the
following BSDE:
\bel{BSDE1}\left\{\2n\ba{ll}
\ns\ds dY_1(s)=-\big[A(s)^\top Y_1(s)+C(s)^\top
Z_1(s)\big]ds+Z_1(s)dW(s),\qq s\in[t,T],\\
\ns\ds Y_1(T)=\eta\in\cX_T.\ea\right.\ee
Then with the well-defined adjoint operators, we can rewrite the
cost functional as follows:
\bel{J-rep}\ba{ll}
\ns\ds J(t,x;u(\cd))=\lan\big(\h L_t^*G\h
L_t+L_t^*QL_t+SL_t+L_t^*S^\top\2n+R\big)u,u\ran\\
\ns\ds\qq\qq\qq~~+2\lan\big(\h
L_t^*G\h\G_t+L_t^*Q\G_t+S\G_t\big)x,u\ran
+\lan\big(\h\G_t^*G\h\G_t+\G_t^*Q\G_t\big)x,x\ran\\
\ns\ds\qq\qq\qq~~+2\lan x,\h\G_t^*\big[Gh_t(T)+g\big]+\G_t^*\big(Qh_t+q\big)\ran\\
\ns\ds\qq\qq\qq~~+2\lan u,\h L_t^*\big[Gh_t(T)+g\big]+L_t^*\big(Qh_t+q\big)+Sh_t+\rho\ran\\
\ns\ds\qq\qq\qq~~+\lan h_t(T),Gh_t(T)+2g\ran+\lan h_t,Qh_t+2q\ran\\
\ns\ds\qq\qq\q\,\equiv\lan M_2(t)u,u\ran+2\lan M_1(t)x,u\ran+\lan M_0(t)x,x\ran
+2\lan x,y_t\ran+2\lan u,\nu_t\ran+c_t,\ea\ee
with $M_2(t):\cU[t,T]\to\cU[t,T]$ being bounded and self-adjoint,
$M_1(t):\dbR^n\to\cU[t,T]$ being bounded, and $M_0(t)\in\dbS^n$;
$y_t\in\dbR^n$, $\nu_t(\cd)\in\cU[t,T]$ and $c_t\in\dbR$. Note that
$\nu_t(\cd), y_t, c_t=0$ when $b(\cd), \si(\cd), g, q(\cd), \rho(\cd)=0$.
This gives (\ref{J-rep1}). Further, if we let
$$X_0(\cd)=(L_tu)(\cd),\qq Y(\cd)=Y_0(\cd)+Y_1(\cd),\qq Z(\cd)=Z_0(\cd)+Z_1(\cd),$$
with
$$\xi(\cd)=Q(\cd)X_0(\cd)+S(\cd)^\top u(\cd),\qq\eta=GX_0(T),$$
then $X_0(\cd)$ satisfies
\bel{X_0}\left\{\2n\ba{ll}
\ns\ds dX_0(s)=\big[A(s)X_0(s)+B(s)u(s)\big]ds+\big[C(s)X_0(s)+D(s)u(s)\big]dW(s),\qq s\in[t,T],\\
\ns\ds X_0(t)=0,\ea\right.\ee
and $(Y(\cd),Z(\cd))$ is the adapted solution to the following BSDE:
$$\left\{\2n\ba{ll}
\ns\ds dY(s)=-\big[A(s)^\top Y(s)+C(s)^\top
Z(s)+Q(s)X_0(s)+S(s)^\top u(s)\big]ds+Z(s)dW(s),\qq
s\in[t,T],\\
\ns\ds Y(T)=GX_0(T).\ea\right.$$
Thus (\ref{L_2-rep}) follows. If we let
$$\bar X_0(\cd)=(\G_tx)(\cd),\qq\bar Y(\cd)=Y_0(\cd)+Z_1(\cd),\qq
\bar Z(\cd)=Z_0(\cd)+Z_1(\cd),$$
with
$$\xi(\cd)=Q(\cd)\bar X_0(\cd),\qq\eta=G\bar X_0(T),$$
then $\bar X_0(\cd)$ satisfies
\bel{bar-X_0}\left\{\2n\ba{ll}
\ns\ds d\bar X_0(s)=A(s)\bar X_0(s)ds+C(s)\bar X_0(s)dW(s),\qq s\in[t,T],\\
\ns\ds\bar X_0(t)=x,\ea\right.\ee
and $(\bar Y(\cd),\bar Z(\cd))$ is the adapted solution to the following BSDE:
$$\left\{\2n\ba{ll}
\ns\ds d\bar Y(s)=-\big[A(s)^\top\bar Y(s)+C(s)^\top \bar
Z(s)+Q(s)\bar X_0(s)\big]ds+\bar Z(s)dW(s),\qq
s\in[t,T],\\
\ns\ds\bar Y(T)=G\bar X_0(T).\ea\right.$$
Thus (\ref{L_1-rep}) follows. Likewise, if we let
$$\h X_0(\cd)=h_t(\cd),\qq\h Y(\cd)=Y_0(\cd)+Z_1(\cd),\qq\h Z(\cd)=Z_0(\cd)+Z_1(\cd),$$
with
$$\xi(\cd)=Q(\cd)\h X_0(\cd)+q(\cd),\qq\eta=G\h X_0(T)+g,$$
then $\h X_0(\cd)$ satisfies
\bel{h-X_0}\left\{\2n\ba{ll}
\ns\ds d\h X_0(s)=\big[A(s)\h X_0(s)+b(s)\big]ds+\big[C(s)\h X_0(s)+\si(s)\big]dW(s),\qq s\in[t,T],\\
\ns\ds\h X_0(t)=0,\ea\right.\ee
and $(\h Y(\cd),\h Z(\cd))$ is the adapted solution to the following BSDE:
$$\left\{\2n\ba{ll}
\ns\ds d\h Y(s)=-\big[A(s)^\top\h Y(s)+C(s)^\top\h Z(s)+Q(s)\h X_0(s)+q(s)\big]ds+\h Z(s)dW(s),\qq
s\in[t,T],\\
\ns\ds\h Y(T)=G\h X_0(T)+g.\ea\right.$$
Thus (\ref{nu-rep}) follows. Finally, we know that
$$\ba{ll}
\ns\ds\lan M_0(t)x,x\ran=J^0(t,x;0)=\dbE\left[\lan G\bar X_0(T),\bar
X_0(T)\ran+\int_t^T\lan Q(s)\bar X_0(s),\bar X_0(s)\ran ds\right]\\
\ns\ds=\dbE\left[\lan
G\F(T)\F(t)^{-1}x,\F(T)\F(t)^{-1}x\ran+\int_t^T\lan
Q(s)\F(s)\F(t)^{-1}x,\F(s)\F(t)^{-1}x\ran ds\right].\ea$$
Then (\ref{L_0}) follows. Also, by Lemma 2.4, we see that
$M_0(\cd)$ is the unique solution of Lyapunov equation (\ref{3.8}).
\endpf

\ms

From the representation of the cost functional, we have the
following simple corollary.

\ms

\bf Corollary 3.2. \sl Let {\rm(H1)--(H2)} hold and $t\in[0,T)$ be
given. For any $x\in\dbR^n, \l\in\dbR$ and $u(\cd),
v(\cd)\in\cU[t,T]$, the following holds:
\bel{u+v-1}\ba{ll}
\ns\ds J(t,x;u(\cd)+\l v(\cd))=J(t,x;u(\cd))+\l^2J^0(t,0;v(\cd))\\
\ns\ds\qq\qq\qq\qq\q~+2\l\dbE\int_t^T\lan B(s)^\top Y(s)\1n+\1nD(s)^\top
Z(s)\1n+\1nS(s)X(s)\1n +\1nR(s)u(s)\1n+\1n\rho(s),v(s)\ran ds,\ea\ee
where $(X(\cd),Y(\cd),Z(\cd))$ is the adapted solution to the
following (decoupled) linear FBSDE:
\bel{FBSDE03}\left\{\2n\ba{ll}
\ns\ds
dX(s)=\big[A(s)X(s)+B(s)u(s)+b(s)\big]ds\\
\ns\ds\qq\qq~+\big[C(s)X(s)+D(s)u(s)+\si(s)\big]dW(s),\qq
s\in[t,T],\\
\ns\ds dY(s)=-\big[A(s)^\top Y(s)+C(s)^\top Z(s)+Q(s)X(s)+S(s)^\top u(s)+q(s)\big]ds\\
\ns\ds\qq\qq\qq\qq\qq\qq\qq\qq\q~+Z(s)dW(s),\qq s\in[t,T], \\
\ns\ds X(t)=x,\qq Y(T)=GX(T)+g.\ea\right.\ee
Consequently, the map $u(\cd)\mapsto J(t,x;u(\cd))$ is Fr\'echet
differentiable with the Fr\'echet derivative given by
\bel{DJ}\cD J(t,x;u(\cd))(s)=2\big[B(s)^\top Y(s)+D(s)^\top
Z(s)+S(s)X(s)+R(s)u(s)+\rho(s)\big],\qq s\in[t,T],\ee
and $(\ref{u+v-1})$ can also be written as
\bel{u+v-1*}J(t,x;u(\cd)+\l
v(\cd))=J(t,x;u(\cd))+\l^2J^0(t,0;v(\cd))+\l\dbE\int_t^T\lan\cD
J(t,x;u(\cd))(s),v(s)\ran ds.\ee

\ms

\it Proof. \rm From Proposition 3.1, we have
$$\ba{ll}
\ns\ds J(t,x;u(\cd)+\l v(\cd))\\
\ns\ds=\lan M_2(t)(u+\l v),u+\l v\ran+2\lan M_1(t)x,u+\l v\ran+\lan M_0(t)x,x\ran
+2\lan u+\l v,\nu_t\ran+2\lan x,y_t\ran+c_t\\
\ns\ds=\lan M_2(t)u,u\ran+2\l\lan M_2(t)u,v\ran+\l^2\lan
M_2(t)v,v\ran+2\lan M_1(t)x,u\ran
+2\l\lan M_1(t)x,v\ran+\lan M_0(t)x,x\ran\\
\ns\ds\q~+2\lan u,\nu_t\ran+2\l\lan v,\nu_t\ran+2\lan x,y_t\ran+c_t\\
\ns\ds=J(t,x;u(\cd))+\l^2J^0(t,0;v(\cd))+2\l\lan M_2(t)u+M_1(t)x+\nu_t,v\ran.\ea$$
From the representation of $M_1(t)$, $M_2(t)$ and $\nu_t$ in Proposition 3.1,
we see that
$$(M_2(t)u)(s)+(M_1(t)x)(s)+\nu_t(s)=B(s)^\top Y(s)+D(s)^\top Z(s)+S(s)X(s)+R(s)u(s)+\rho(s),\q s\in[t,T],$$
with $(X(\cd),Y(\cd),Z(\cd))$ being the adapted solution to the
FBSDE (\ref{FBSDE03}). The rest of the proof is clear. \endpf

\ms

Note that if $u(\cd)$ happens to be an open-loop optimal control of
Problem (SLQ), then the following {\it stationarity condition}
holds:
\bel{J_u=0}\cD J(t,x;u(\cd))=2\big[B(s)^\top Y(s)+D(s)^\top
Z(s)+S(s)X(s)+R(s)u(s)+\rho(s)\big]=0,\qq s\in[t,T],\ee
which brings a coupling into the FBSDE (\ref{FBSDE03}). We call
(\ref{FBSDE03}), together with the stationarity condition
(\ref{J_u=0}), the {\it optimality system} for the open-loop optimal
control of Problem (SLQ).

\ms

The following is concerned with the convexity of the cost
functional, whose proof is straightforward, by making use of the
representation (\ref{J-rep1}) of the cost functional.

\ms

\bf Corollary 3.3. \sl Let {\rm(H1)--(H2)} hold and let $t\in[0,T)$
be given. Then the following are equivalent:

\ms

{\rm(i)} $u(\cd)\mapsto J(t,x;u(\cd))$ is convex, for some
$x\in\dbR^n$.

\ms

{\rm(ii)} $u(\cd)\mapsto J(t,x;u(\cd))$ is convex, for any
$x\in\dbR^n$.

\ms

{\rm(iii)} $u(\cd)\mapsto J^0(t,x;u(\cd))$ is convex, for some
$x\in\dbR^n$.

\ms

{\rm(iv)} $u(\cd)\mapsto J^0(t,x;u(\cd))$ is convex, for any
$x\in\dbR^n$.

\ms

{\rm(v)} $J^0(t,0;u(\cd))\ges0$, for all $u(\cd)\in\cU[t,T]$.

\ms

{\rm(vi)} $M_2(t)\ges0$.

\rm

\ms

Similar to the above, we have that $u(\cd)\mapsto J(t,x;u(\cd))$ is
uniformly convex if and only if
\bel{J>l}J^0(t,0;u(\cd))\ges\l\dbE\int_t^T|u(s)|^2ds,\qq\forall
u(\cd)\in\cU[t,T],\ee
for some $\l>0$. This is also equivalent to the following:
\bel{M_2>l}M_2(t)\ges\l I,\ee
for some $\l>0$. Further, it is obvious that if the standard
conditions (\ref{classical}) hold, then
\bel{}M_2(t)=\h L_t^*G\h L_t+L_t^*(Q-S^\top R^{-1}S)L_t+(L_t^*S^\top
R^{-{1\over2}}+R^{1\over2})(R^{-{1\over2}}SL_t+R^{1\over2})\ges0,\ee
which means that the functional $u(\cd)\mapsto J^0(t,0,u(\cd))$ is
convex. The following result tells us that under (\ref{classical}),
one actually has the uniform convexity of the cost functional.

\ms

\bf Proposition 3.4. \sl Let {\rm(H1)--(H2)} and {\rm(\ref{classical})} hold. Then for any $t\in[0,T)$, the
map $u(\cd)\mapsto J^0(t,0;u(\cd))$ is uniformly convex.

\ms

\it Proof. \rm For any $u(\cd)\in\cU[t,T]$, let $X^{(u)}(\cd)$ be the solution of
$$\left\{\2n\ba{ll}
\ns\ds dX^{(u)}(s)=\big[A(s)X^{(u)}(s)+B(s)u(s)\big]ds+\big[C(s)X^{(u)}(s)+D(s)u(s)\big]dW(s),\qq s\in[t,T], \\
\ns\ds X^{(u)}(t)=0.\ea\right.$$
Then by Lemma 2.5 (taking $\Th(\cd)=-R(\cd)^{-1}S(\cd)$), we have
$$\ba{ll}
\ns\ds J^0(t,0;u(\cd))=\dbE\left\{\langle GX^{(u)}(T),X^{(u)}(T)\rangle
+\int_t^T\[\langle QX^{(u)},X^{(u)}\rangle+2\langle SX^{(u)},u\rangle+\langle Ru,u\rangle\]ds\right\}\\
\ns\ds\qq\qq\q~\ges\dbE\int_t^T\[\langle QX^{(u)},X^{(u)}\rangle+2\langle SX^{(u)},u\rangle+\langle Ru,u\rangle\]ds\\
\ns\ds\qq\qq\q~=\dbE\int_t^T\[\lan\big(Q\1n-\1nS^\top R^{-1}S\big)X^{(u)},X^{(u)}\ran\1n+\1n\lan R\big(u\1n+\1nR^{-1}SX^{(u)}\big),u\1n+\1nR^{-1}SX^{(u)}\ran\]ds\\
\ns\ds\qq\qq\q~\ges\d\dbE\int_t^T\big|u+R^{-1}SX^{(u)}\big|^2
ds\ges\d\g\dbE\int_t^T|u(s)|^2 ds,\qq \forall
u(\cd)\in\cU[t,T],\ea$$
for some $\g>0$. This completes the proof.
\endpf

\section{Solvabilities of Problem (SLQ), Uniform Convexity of the Cost Functional, and the Riccati Equation}

We begin with a simple result concerning the open-loop solvability
of Problem (SLQ).

\ms

\bf Proposition 4.1. \sl Let {\rm(H1)--(H2)} hold. Suppose the map
$u(\cd)\mapsto J^0(0,0;u(\cd))$ is uniformly convex. Then Problem
{\rm(SLQ)} is uniquely open-loop solvable, and there exists a
constant $\a\in\dbR$ such that
\bel{uni-convex-prop0}V^0(t,x)\ges\a|x|^2,\qq\forall
(t,x)\in[0,T]\times\dbR^n.\ee

\ms

\rm

Note that in the above, the constant $\a$ does not have to be
nonnegative.

\ms

\it Proof. \rm First of all, by the uniform convexity of
$u(\cd)\mapsto J^0(0,0;u(\cd))$, we may assume that
\bel{J>l*}J^0(0,0;u(\cd))\ges\l\,\dbE\1n\int_0^T|u(s)|^2ds,\qq\forall
u(\cd)\in\cU[0,T],\ee
for some $\l>0$. Now, for any $t\in[0,T)$, and any
$u(\cd)\in\cU[t,T]$, we define the {\it zero-extension} of $u(\cd)$
as follows:
\bel{ext}[\,0I_{[0,t)}\oplus u(\cd)](s)=\left\{\2n\ba{ll}0,\qq\ s\in[0,t),\\
\ns\ds u(s),\q s\in[t,T].\ea\right.\ee
Then $v(\cd)\equiv0I_{[0,t)}\oplus u(\cd)\in\cU[0,T]$, and due to
the initial state being 0, the solution $X(s)$ of
$$\left\{\2n\ba{ll}
\ns\ds dX(s)=\big[A(s)X(s)+B(s)v(s)\big]ds+\big[C(s)X(s)+D(s)v(s)\big]dW(s),\qq s\in[0,T], \\
\ns\ds X(0)=0,\ea\right.$$
satisfies
$$X(s)=0,\qq s\in[0,t].$$
Hence,
\bel{4.4}J^0(t,0;u(\cd))=J^0(0,0;0I_{[0,t)}\oplus u(\cd))\ges
\l\,\dbE\1n\int_0^T\big|[0I_{[0,t)}\oplus
u(\cd)](s)\big|^2ds=\l\,\dbE\1n\int_t^T|u(s)|^2ds.\ee
Thus, $u(\cd)\mapsto J^0(t,x;u(\cd))$ is uniformly convex for any
given $(t,x)\in[0,T)\times\dbR^n$. By Corollary 3.2, we have
\bel{uni-convex-prop1}\ba{ll}
\ns\ds J(t,x;u(\cd))=J(t,x;0)+J^0(t,0;u(\cd))+\dbE\int_t^T\lan\cD J(t,x;0)(s),u(s)\ran ds\\
\ns\ds\ges J(t,x;0)+J^0(t,0;u(\cd))-{\l\over2}\dbE\int_t^T|u(s)|^2ds-{1\over2\l}\dbE\int_t^T|\cD J(t,x;0)(s)|^2ds\\
\ns\ds\ges{\l\over2}\dbE\int_t^T|u(s)|^2ds+J(t,x;0)-{1\over2\l}\dbE\int_t^T|\cD
J(t,x;0)(s)|^2ds,\qq\forall u(\cd)\in\cU[t,T].\ea\ee
Consequently, by a standard argument involving minimizing sequence
and locally weak compactness of Hilbert spaces, we see that for any
given initial pair $(t,x)\in [0,T)\times\dbR^n$, Problem (SLQ)
admits a unique open-loop optimal control. Moreover, when $b(\cd), \si(\cd), g, q(\cd), \rho(\cd)=0$,
(\ref{uni-convex-prop1}) implies that
\bel{uni-convex-prop2}V^0(t,x)\ges J^0(t,x;0)-{1\over2\l}\dbE\int_t^T|\cD J^0(t,x;0)(s)|^2ds.\ee
Note that the functions on the right-hand side of (\ref{uni-convex-prop2}) are
quadratic in $x$ and continuous in $t$. (\ref{uni-convex-prop0})
follows immediately. \endpf

\ms

Now, we introduce the following Riccati equation associated with
Problem (SLQ):
\bel{Riccati}\left\{\2n\ba{ll}
\ns\ds\dot P(s)+P(s)A(s)+A(s)^\top P(s)+C(s)^\top P(s)C(s)+Q(s)\\
\ns\ds\q-\big[P(s)B(s)+C(s)^\top P(s)D(s)+S(s)^\top\big]\big[R(s)+D(s)^\top P(s)D(s)\big]^\dag\\
\ns\ds\qq\cd\big[B(s)^\top P(s)+D(s)^\top P(s)C(s)+S(s)\big]=0,\qq \ae~s\in[0,T],\\
\ns\ds P(T)=G.\ea\right.\ee
A solution $P(\cd)\in C([0,T];\dbS^n)$ of (\ref{Riccati}) is said to
be {\it regular} if
\bel{regular-1}\cR\big(B(s)^\top P(s)+D(s)^\top
P(s)C(s)+S(s)\big)\subseteq\cR\big(R(s)+D(s)^\top P(s)D(s)\big),\qq
\ae~s\in[0,T],\ee
\bel{regular-2}\big[R(\cd)+D(\cd)^\top
P(\cd)D(\cd)\big]^\dag\big[B(\cd)^\top P(\cd) +D(\cd)^\top
P(\cd)C(\cd)+S(\cd)\big]\in L^2(0,T;\dbR^{m\times n}),\ee
and
\bel{regular-3}R(s)+D(s)^\top P(s)D(s)\ges0,\qq\ae~s\in[0,T].\ee
A solution $P(\cd)$ of (\ref{Riccati}) is said to be {\it strongly
regular} if
\bel{strong-regular}R(s)+D(s)^\top P(s)D(s)\ges \l  I,\qq\ae~s\in[0,T],\ee
for some $\l>0$. The Riccati equation (\ref{Riccati}) is said to be
({\it strongly}) {\it regularly solvable}, if it admits a (strongly)
regular solution. Clearly, condition (\ref{strong-regular}) implies
(\ref{regular-1})--(\ref{regular-3}). Thus, a strongly regular
solution $P(\cd)$ must be regular. Moreover, it was shown in
\cite{Sun-Yong 2014} that if a regular solution of (\ref{Riccati})
exists, it must be unique.

\ms

In \cite{Ait Rami-Moore-Zhou 2001}, it was showed that for Problem
{$\rm(SLQ)^0$}, the existence of a continuous open-loop optimal
control for any initial pair $(t,x)\in[0,T]\times\dbR^n$ is
equivalent to the solvability of the corresponding Riccati equation
{\rm(\ref{Riccati})} with constraints (\ref{regular-1}) and
(\ref{regular-3}). More precisely, their result can be stated as
follows (in terms of our notations and equation numbers):

\ms

``\it\textbf{Theorem} {\rm 4.2}. \it Suppose that $B(\cd), C(\cd), D(\cd), Q(\cd), R(\cd)$ are continuous and $S(\cd)=0$.
Then Problem {$\rm(SLQ)^0$} has a continuous open-loop optimal control for any initial pair
$(t,x)\in[0,T]\times\dbR^n$ if and only if the Riccati equation {\rm(\ref{Riccati})} has a solution $P(\cd)$
such that {\rm(\ref{regular-1})} and {\rm(\ref{regular-3})} hold.\rm"

\ms

\no This result is incorrect. Here is a simple counterexample.

\ms

\bf Example 4.3. \rm Consider the following deterministic one-dimensional state equation:
$$\left\{\2n\ba{ll}
\ns\ds dX(s)=u(s)ds,\qq s\in[t,1], \\
\ns\ds X(t)=x,\ea\right.$$
and cost functional
$$J(t,x;u(\cd))=X(1)^2+\int_t^1s^2u(s)^2ds.$$
The Riccati equation for the above problem reads
\bel{ex4.2-Ric}\left\{\2n\ba{ll}
\ns\ds\dot P(t)={P(t)^2\over t^2}{\bf 1}_{(0,1]}(t),\qq\ae~t\in[0,1]\\
\ns\ds P(1)=1.\ea\right.\ee
It is easy to see that $P(t)=t$ is the unique solution of (\ref{ex4.2-Ric}),
satisfying {\rm(\ref{regular-1})} and {\rm(\ref{regular-3})}. Now,
we claim that this problem does not admit an open-loop optimal
control for initial pair $(0,x)$ with $x\neq0$. In fact, if for some
$x\neq0$, there exists an open-loop optimal control
$u^*(\cd)\in\cU[0,T]$, then by the maximum principle, the solution
$(X^*(\cd),Y^*(\cd))$ of the following (decoupled) forward-backward
differential equation:
\bel{ex4.2-FBDE}\left\{\2n\ba{ll}
\ns\ds \dot X^*(s)=u^*(s),\q\dot Y^*(s)=0,\qq s\in[0,1],\\
\ns\ds X^*(0)=x,\q Y^*(1)=X^*(1),\ea\right.\ee
must satisfy
\bel{ex4.2-constraint}Y^*(s)+s^2u^*(s)=0,\qq\ae~s\in[0,1].\ee
Observe that the solution $(X^*(\cd),Y^*(\cd))$ of (\ref{ex4.2-FBDE}) is given by
$$X^*(s)=x+\int_0^su^*(r)dr,\q Y^*(s)\equiv X^*(1),\qq s\in[0,1].$$
Hence,
$$u^*(s)={X^*(1)\over s^2},\qq\ae~s\in(0,1].$$
Noting that $u^*(\cd)$ is square-integrable, we must have $X^*(1)=0$ and hence $u^*(\cd)=0$. Consequently,
$$x=X^*(1)-\int_0^1u^*(r)dr=0,$$
which is a contradiction.

\ms

From the above example, we see that the sufficiency part of the
above Theorem 4.2 (a result from \cite{Ait Rami-Moore-Zhou 2001})
does not hold. We will see in Section 7 that the necessity part does
not hold either.

\ms

Instead of Theorem 4.2, in \cite{Sun-Yong 2014}, the following were
proved, which establishes the equivalence between the closed-loop
solvability of Problem (SLQ) and the regular solvability of the
Riccati equation (\ref{Riccati}).

\ms


\bf Theorem 4.4. \sl Let {\rm(H1)--(H2)} hold. Problem {\rm(SLQ)} is
closed-loop solvable on $[0,T]$ if and only if the
Riccati equation {\rm(\ref{Riccati})} admits a regular solution
$P(\cd)\in C([0,T];\dbS^n)$ and the adapted solution $(\eta(\cd),\z(\cd))$ of the
following BSDE:
\bel{eta-zeta}\left\{\2n\ba{ll}
\ns\ds d\eta(s)=-\Big\{\big[A^\top\2n-(PB+C^\top\1n PD+S^\top)(R+D^\top\1n
PD)^\dag B^\top\big]\eta\\
\ns\ds\qq\qq\q+\big[C^\top\2n-(PB+C^\top\1n PD+S^\top)(R+D^\top\1n PD)^\dag D^\top\big]\z\\
\ns\ds\qq\qq\q+\big[C^\top\2n-(PB+C^\top\1n PD+S^\top)(R+D^\top\1n PD)^\dag D^\top\big]P\si\\
\ns\ds\qq\qq\q-(PB+C^\top\1n PD+S^\top)(R+D^\top\1n PD)^\dag
\rho+Pb+q\Big\}ds+\z dW(s),\q s\in[0,T],\\
\ns\ds\eta(T)=g,\ea\right.\ee
satisfies
\bel{eta-zeta-regularity}\left\{\2n\ba{ll}
\ns\ds B^\top\eta+D^\top\z+D^\top P\si+\rho\in\cR(R+D^\top PD), \qq \ae~\as\\
\ns\ds (R+D^\top PD)^\dag(B^\top\eta+D^\top\z+D^\top P\si+\rho)\in L_\dbF^2(0,T;\dbR^m).\ea\right.\ee
In this case, Problem {\rm(SLQ)} is
closed-loop solvable on any $[t,T]$, and the closed-loop optimal
strategy $(\Th^*(\cd),v^*(\cd))$ admits the following
representation:
\bel{Th-v-rep}\left\{\2n\ba{ll}
\Th^*=-(R+D^\top PD)^\dag(B^\top P+D^\top PC+S)
+\big[I-(R+D^\top PD)^\dag(R+D^\top PD)\big]\Pi,\\
\ns\ds v^*=-(R+D^\top PD)^\dag(B^\top\eta+D^\top\z+D^\top P\si+\rho)
+\big[I-(R+D^\top PD)^\dag(R+D^\top PD)\big]\n,\ea\right.\ee
for some $\Pi(\cd)\in L^2(t,T;\dbR^{m\times n})$ and $\n(\cd)\in
L_\dbF^2(t,T;\dbR^m)$, and the value function is given by
\bel{Value}\ba{ll}
\ns\ds V(t,x)=\dbE\bigg\{\langle P(t)x,x\rangle+2\langle\eta(t),x\rangle+\int_t^T\[\langle
P\si,\si\rangle+2\langle\eta,b\rangle+2\langle\z,\si\rangle\\
\ns\ds\qq\qq\q\ ~-\lan(R+D^\top PD)^\dag(B^\top\eta+D^\top\z+D^\top P\si+\rho),
B^\top\eta+D^\top\z+D^\top P\si+\rho\ran\]ds\bigg\}.\ea\ee

\rm

Note that in Example 4.3, the solution $P(t)=t$ to the Riccati
equation (\ref{ex4.2-Ric}) is not regular since
$$[R(t)+D(t)^\top P(t)D(t)]^\dag[B(t)^TP(t)+D(t)^\top
P(t)C(t)+S(t)]={1\over t}{\bf 1}_{(0,1]}(t)\notin L^2(0,1;\dbR).$$
Hence, by Theorem 4.4, the corresponding LQ problem does not have a
closed-loop optimal strategy.

\ms

From the above theorem, we see that the existence of a strongly
regular solution to the Riccati equation (\ref{Riccati}) implies the
unique closed-loop solvability of Problem (SLQ), which, by the
remark right after Definition 2.1, implies the unique open-loop
solvability of Problem (SLQ). Particularly, when $b(\cd), \si(\cd),
g, q(\cd), \rho(\cd)=0$, the adapted solution of (\ref{eta-zeta}) is
$(0,0)$, and (\ref{eta-zeta-regularity}) holds automatically. Thus,
the existence of a regular solution to the Riccati equation
(\ref{Riccati}) is equivalent to the closed-loop solvability of
Problem {$\rm(SLQ)^0$}, which implies the open-loop solvability of
Problem {$\rm(SLQ)^0$}. However, from Example 2.2, we know that the
inverse is false, i.e., it may happen that for any initial pair
$(t,x)\in[0,T)\times\dbR^n$, Problem {$\rm(SLQ)^0$} admits an
open-loop optimal control (which could even be continuous), while
the problem is not closed-loop solvable, which means that
(\ref{Riccati}) does not admit a regular solution (see Section 7 for
further details). On the other hand, it is known that under the
standard conditions (\ref{classical}), the Riccati equation
(\ref{Riccati}) admits a unique solution $P(\cd)\in
C([0,T];\cl{\dbS_+^n})$, and Problem (SLQ) admits a unique open-loop
optimal control which has a state feedback form, represented via the
solution of the Riccati equation (see \cite{Yong-Zhou 1999,
Chen-Zhou 2000}). To summarize up, we have the following diagram:

{\footnotesize
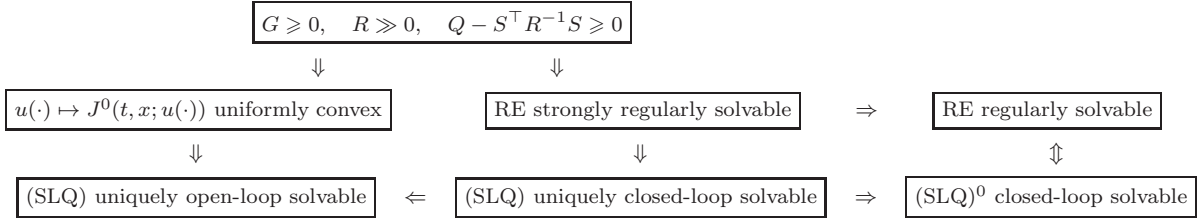
$$\ba{ll}
\ns\ds\qq\qq\qq\qq\qq\q\boxed{G\ges0,\q R\gg0,\q Q-S^\top R^{-1}S\ges0}\\
\ns\ds\qq\qq\qq\qq\qq\qq\q~\Downarrow\qq\qq\qq\qq\qq\Downarrow\\
\ns\ds \boxed{\hb{$u(\cd)\mapsto J^0(t,x;u(\cd))$ uniformly
convex}}\qq\qq\,\boxed{\sc\hb{RE strongly regularly solvable}}
\qq\,\Rightarrow\qq\,\boxed{\sc\hb{RE regularly solvable}}\\
\ns\ds\qq\qq\qq\qq\,\Downarrow\qq\qq\qq\qq\qq\qq\qq\qq\qq\q\,\Downarrow\qq\qq\qq\qq\qq\qq\qq\qq\q~\Updownarrow\\
\ns\ds~\,\boxed{\hb{(SLQ) uniquely open-loop
solvable}}\q\Leftarrow\q\boxed{\hb{(SLQ) uniquely closed-loop
solvable}}\q\Rightarrow\q\boxed{\hb{$\rm(SLQ)^0$ closed-loop
solvable}}\ea$$}

\no where ``RE" stands for the Riccati equation (\ref{Riccati}). It
is clear that the uniform convexity of the map $u(\cd)\mapsto
J^0(t,x;u(\cd))$ does not imply the standard conditions
(\ref{classical}), which will be even clearer by the results of
Section 5 below. Therefore, it is a desire to establish the
following:
$$\boxed{\hb{$u(\cd)\mapsto J^0(t,x;u(\cd))$ uniformly
convex}}\q\iff\q\boxed{\hb{RE strongly regularly solvable}}$$
This is our next goal. To achieve this, we first present the
following proposition, which plays a key technical role in this
section.

\ms

\bf Proposition 4.5. \sl Let {\rm(H1)--(H2)} and {\rm(\ref{J>l*})}
hold. Then for any $\Th(\cd)\in L^2(0,T;\dbR^{m\times n})$, the
solution $P(\cd)\in C([0,T];\dbS^n)$ to the following Lyapunov equation:
\bel{P-Th}\left\{\2n\ba{ll}
\ns\ds\dot P+P(A+B\Th)+(A+B\Th)^\top P+(C+D\Th)^\top P(C+D\Th)\\
\ns\ds\q+\,\Th^\top R\Th+S^\top\Th+\Th^\top S+Q=0,\qq\ae~s\in[0,T],\\
\ns\ds P(T)=G,\ea\right.\ee
satisfies
\bel{Convex-prop-1}R(t)+D(t)^\top P(t)D(t)\ges\l I, \q\ae~t\in[0,T],\qq\hb{and}\qq P(t)\ges \a I,\q\forall t\in[0,T],\ee
where $\a\in\dbR$ is the constant appears in
{\rm(\ref{uni-convex-prop0})}.

\ms

\it Proof. \rm Let $\Th(\cd)\in L^2(0,T;\dbR^{m\times n})$ and let
$P(\cd)$ be the solution to {\rm(\ref{P-Th})}. For any
$u(\cd)\in\cU[0,T]$, let $X_0(\cd)$ be the solution of
$$\left\{\2n\ba{ll}
\ns\ds dX_0(s)=\big[(A+B\Th) X_0+Bu\big]ds+\big[(C+D\Th)X_0+Du\big]dW(s),\qq s\in[0,T], \\
\ns\ds X_0(0)=0.\ea\right.$$
By {\rm(\ref{J>l*})} and Lemma 2.3, we have
$$\ba{ll}
\ns\ds \l\dbE\int_0^T|\Th(s)X_0(s)+u(s)|^2ds\les J^0(0,0;\Th(\cd)X_0(\cd)+u(\cd))\\
\ns\ds=\dbE\int_0^T\Big\{2\lan\big[B^\top P+D^\top PC+S+(R+D^\top PD)\Th\big]X_0,u\ran+\lan (R+D^\top PD)u,u\ran\Big\}ds.\ea$$
Hence, for any $u(\cd)\in\cU[0,T]$, the following holds:
\bel{P>LI}\ba{ll}
\ns\ds\dbE\int_0^T\Big\{2\lan\big[B^\top P+D^\top PC+S+(R+D^\top PD-\l I)\Th\big]X_0,u\ran\\
\ns\ds\qq\q~+\lan (R+D^\top PD-\l I)u,u\ran\Big\}ds=\l\dbE\int_0^T|\Th(s)X_0(s)|^2ds\ges0.\ea\ee
Now, fix any $u_0\in\dbR^m$, take $u(s)=u_0{\bf 1}_{[t,t+h]}(s)$,
with $0\les t<t+h\les T$. Then
$$\left\{\2n\ba{ll}
\ns\ds d\big[\dbE X_0(s)\big]=\Big\{\big[A(s)+B(s)\Th(s)\big]\dbE X_0(s)+B(s)u_0{\bf1}_{[t,t+h]}(s)\Big\}ds,\qq
s\in[0,T],\\
\ns\ds\dbE X_0(0)=0.\ea\right.$$
Hence,
$$\dbE X_0(s)=\left\{\2n\ba{ll}0,\qq\qq\qq\qq\qq\qq\qq\q\1n s\in[0,t],\\
\ns\ds\F(s)\int_t^{s\land(t+h)}\F(r)^{-1}B(r)u_0dr,\qq
s\in[t,T],\ea\right.$$
where $\F(\cd)$ is the solution of the following $\dbR^{n\times n}$-valued ordinary differential equation:
$$\left\{\2n\ba{ll}
\ns\ds \dot\F(s)=\big[A(s)+B(s)\Th(s)\big]\F(s),\qq s\in[0,T], \\
\ns\ds\F(0)=I.\ea\right.$$
Consequently, (\ref{P>LI}) becomes
$$\ba{ll}
\ns\ds\int_t^{t+h}\Big\{2\lan\big[B^\top P+D^\top PC+S
+(R+D^\top PD-\l I)\Th\big]\F(s)\int_t^s\F(r)^{-1}B(r)u_0dr,u_0\ran\\
\ns\ds\qq\q~+\lan (R+D^\top PD-\l I)u_0,u_0\ran\Big\}ds\ges0.\ea$$
Dividing both sides of the above by $h$ and letting $h\to 0$, we obtain
$$\lan\big[R(t)+D(t)^\top P(t)D(t)-\l I\big]u_0,u_0\ran\ges 0,\qq\ae~t\in[0,T],\q \forall u_0\in\dbR^m.$$
The first inequality in (\ref{Convex-prop-1}) follows. To prove the second, for any $(t,x)\in[0,T)\times\dbR^n$ and $u(\cd)\in\cU[t,T]$,
let $X(\cd)$ be the solution of
$$\left\{\2n\ba{ll}
\ns\ds dX(s)=\big[(A+B\Th) X+Bu\big]ds+\big[(C+D\Th)X+Du\big]dW(s),\qq s\in[t,T], \\
\ns\ds X(t)=x.\ea\right.$$
By Proposition 4.1 and Lemma 2.3, we have
$$\ba{ll}
\ns\ds \a|x|^2\les V^0(t,x)\les J^0(t,x;\Th(\cd)X(\cd)+u(\cd))\\
\ns\ds\qq~\1n=\langle P(t)x,x\rangle\1n+\dbE\int_t^T\Big\{
2\lan\big[B^\top\1nP\1n+\1nD^\top\1n PC\1n+\1nS\1n+\1n(R\1n+\1nD^\top\1n PD)\Th\big]X,u\ran\1n
+\1n\lan(R\1n+\1nD^\top\1n PD)u,u\ran\Big\}ds.\ea$$
In particular, by taking $u(\cd)=0$ in the above, we obtain
$$\langle P(t)x,x\rangle\ges\a|x|^2,\qq\forall (t,x)\in[0,T]\times\dbR^n,$$
and the second inequality therefore follows. \endpf

\ms

Now we state the main result of this section.

\ms

\bf Theorem 4.6. \sl Let {\rm(H1)--(H2)} hold. Then the following statements are equivalent:

\ms

{\rm(i)} The map $u(\cd)\mapsto J^0(0,0;u(\cd))$ is uniformly convex, i.e., there exists a
$\l>0$ such that {\rm(\ref{J>l*})} holds.

\ms

{\rm(ii)} The Riccati equation {\rm(\ref{Riccati})} admits a strongly regular solution $P(\cd)\in C([0,T];\dbS^n)$.

\ms

\it Proof. \rm (i) $\Ra$ (ii). Let $P_0$ be the solution of
\bel{}\left\{\2n\ba{ll}
\ns\ds\dot P_0+P_0A+A^\top P_0+C^\top P_0C+Q=0,\qq\ae~s\in[0,T],\\
\ns\ds P_0(T)=G.\ea\right.\ee
Applying Proposition 4.5 with $\Th=0$, we obtain that
$$R(t)+D(t)^\top P_0(t)D(t)\ges\l I,\q P_0(t)\ges\a I,\qq\ae~t\in[0,T].$$
Next, inductively, for $i = 0,1,2, \cdots$, we set
\bel{Iteration-i}\Th_i=-(R+D^\top P_iD)^{-1}(B^\top P_i+D^\top P_iC+S),\qq A_i=A+B\Th_i,\qq C_i=C+D\Th_i,\ee
and let $P_{i+1}$ be the solution of
\bel{}\left\{\2n\ba{ll}
\ns\ds\dot{P}_{i+1}+P_{i+1}A_i+A_i^\top P_{i+1}+C_i^\top P_{i+1}C_i
+\Th_i^\top R\Th_i+S^\top\Th_i+\Th_i^\top S+Q=0,\qq\ae~s\in[0,T],\\
\ns\ds P_{i+1}(T)=G.\ea\right.\ee
By Proposition 4.5, we see that
\bel{R+Pi-lowerbound}R(t)+D(t)^\top P_{i+1}(t)D(t)\ges\l I,\q
P_{i+1}(t)\ges\a I,\qq \ae~s\in[0,T],\qq i=0,1,2,\cdots.\ee
We now claim that $\{P_i\}_{i=1}^\i$ converges uniformly in
$C([0,T];\dbS^n)$. To show this, let
$$\D_i\deq P_i-P_{i+1},\qq \L_i\deq\Th_{i-1}-\Th_i,\qq i\ges1.$$
Then for $i\ges1$, we have
\bel{Di-equa1}\ba{ll}
\ns\ds -\dot \D_i=\dot{P}_{i+1}-\dot{P}_i \\
\ns\ds\qq~\2n=P_iA_{i-1}+A_{i-1}^\top P_i+C_{i-1}^\top P_iC_{i-1}+\Th_{i-1}^\top R\Th_{i-1}+S^\top\Th_{i-1}+\Th_{i-1}^\top S\\
\ns\ds\qq\q~-P_{i+1}A_i-A_i^\top P_{i+1}-C_i^\top P_{i+1}C_i-\Th_i^\top R\Th_i-S^\top\Th_i-\Th_i^\top S\\
\ns\ds\qq~\2n=\D_iA_i+A_i^\top\D_i+C_i^\top\D_iC_i+P_i(A_{i-1}-A_i)+(A_{i-1}-A_i)^\top P_i\\
\ns\ds\qq\q~+C_{i-1}^\top P_iC_{i-1}-C_i^\top P_iC_i+\Th_{i-1}^\top R\Th_{i-1}-\Th_i^\top R\Th_i+S^\top\L_i+\L_i^\top S.\ea\ee
By (\ref{Iteration-i}), we have the following:
\bel{Di-equa2}\left\{\2n\ba{ll}
\ns\ds A_{i-1}-A_i=B\L_i,\qq C_{i-1}-C_i=D\L_i,\\
\ns\ds C_{i-1}^\top P_iC_{i-1}-C_i^\top P_iC_i=\L_i^\top D^\top P_iD\L_i+C_i^\top P_iD\L_i+\L_i^\top D^\top P_iC_i,\\
\ns\ds \Th_{i-1}^\top R\Th_{i-1}-\Th_i^\top R\Th_i=\L_i^\top R\L_i+\L_i^\top R\Th_i+\Th_i^\top R\L_i.\ea\right.\ee
Note that
$$B^\top P_i+D^\top P_iC_i+R\Th_i+S=B^\top P_i+D^\top P_iC+S+(R+D^\top P_iD)\Th_i=0.$$
Thus, plugging (\ref{Di-equa2}) into (\ref{Di-equa1}) yields
\bel{Di-equa3}\ba{ll}
\ns\ds-\,(\dot\D_i+\Delta_iA_i+A_i^\top\D_i+C_i^\top\D_iC_i)\\
\ns\ds=P_iB\L_i+\L_i^\top B^\top P_i+\L_i^\top D^\top P_iD\L_i+C_i^\top P_iD\L_i+\L_i^\top D^\top P_iC_i\\
\ns\ds\q+\L_i^\top R\L_i+\L_i^\top R\Th_i+\Th_i^\top R\L_i+S^\top\L_i+\L_i^\top S\\
\ns\ds=\L_i^\top(R+D^\top P_iD)\L_i+(P_iB+C_i^\top P_iD+\Th_i^\top R+S^\top)\L_i+\L_i^\top(B^\top P_i+D^\top P_iC_i+R\Th_i+S)\\
\ns\ds=\L_i^\top(R+D^\top P_iD)\L_i\ges0.\ea\ee
Noting that $\D_i(T)=0$ and using Lemma 2.4, also noting (\ref{R+Pi-lowerbound}), we obtain
$$P_1(s)\ges P_i(s)\ges P_{i+1}(s)\ges\a I,\qq\forall s\in [0,T],\q\forall i\ges1.$$
Therefore, the sequence $\{P_i\}_{i=1}^\i$ is uniformly bounded.
Consequently, there exists a constant $K>0$ such that (noting
(\ref{R+Pi-lowerbound}))
\bel{Di-equa4}\left\{\2n\ba{ll}
\ns\ds|P_i(s)|,\ |R_i(s)|\les K,\\
\ns\ds|\Th_i(s)|\les K\big(|B(s)|+|C(s)|+|S(s)|\big),\\
\ns\ds|A_i(s)|\les |A(s)|+K|B(s)|\big(|B(s)|+|C(s)|+|S(s)|\big),\\
\ns\ds|C_i(s)|\les |C(s)|+K\big(|B(s)|+|C(s)|+|S(s)|\big),\ea\right.
\qq\ae~s\in [0,T],\q\forall \ i\ges0,\ee
where $R_i(s)\deq R(s)+D^\top(s)P_i(s)D(s)$. Observe that
\bel{Di-equa5}\ba{ll}
\ns\ds\L_i=\Th_{i-1}-\Th_i \\
\ns\ds\q\ =R_i^{-1}D^\top\D_{i-1}DR_{i-1}^{-1}(B^\top P_i+D^\top P_iC+S)-R_{i-1}^{-1}(B^\top\D_{i-1}+D^\top\D_{i-1}C).\ea\ee
Thus, noting (\ref{Di-equa4}), one has
\bel{3.22}\ba{ll}
\ns\ds|\L_i(s)^\top R_i(s)\L_i(s)|\les\(|\Th_i(s)|+|\Th_{i-1}(s)|\)\,|R_i(s)|\,
|\Th_{i-1}(s)-\Th_i(s)|\\
\ns\ds\qq\qq\qq\qq\2n~\les K\(|B(s)|+|C(s)|+|S(s)|\)^2|\D_{i-1}(s)|.\ea\ee
Equation (\ref{Di-equa3}), together with $\D_i(T)=0$, implies that
$$\D_i(s)=\int^T_s\big(\D_iA_i+A_i^\top\D_i+C_i^\top\D_iC_i+\L_i^\top R_i\L_i\big)dr.$$
Making use of (\ref{3.22}) and still noting (\ref{Di-equa4}), we get
$$|\D_i(s)|\les \int^T_s\f(r)\[|\D_i(r)|+|\D_{i-1}(r)|\]dr,\qq\forall s\in[0,T],\q\forall i\ges1,$$
where $\f(\cd)$ is a nonnegative integrable function independent of $\D_i(\cd)$. By Gronwall's inequality,
$$|\D_i(s)|\les e^{\int_0^T\f(r)dr}\int^T_s\f(r)|\D_{i-1}(r)|dr\equiv c\int^T_s\f(r)|\D_{i-1}(r)|dr.$$
Set
$$a\deq\max_{0\les s\les T}|\D_0(s)|.$$
By induction, we deduce that
$$|\D_i(s)|\les a{c^i\over i!}\(\int_s^T\f(r)dr\)^i,\qq\forall s\in[0,T],$$
which implies the uniform convergence of $\{P_i\}_{i=1}^\i$. We denote $P$ the limit of $\{P_i\}_{i=1}^\infty$, then
(noting (\ref{R+Pi-lowerbound}))
$$R(s)+D(s)^\top P(s)D(s)=\lim_{i\to\i}R(s)+D(s)^\top P_i(s)D(s)\ges\l I,
\qq\ae~s\in[0,T],$$
and as $i\to\infty$,
$$\left\{\2n\ba{ll}
\ns\ds\Th_i\to-(R+D^\top PD)^{-1}(B^\top P+D^\top PC+S)\equiv\Th\q\hb{in $L^2$},\\
\ns\ds A_i\to A+B\Th\q\hb{in $L^1$},\qq C_i\to C+D\Th\q\hb{in $L^2$}.\ea\right.$$
Therefore, $P(\cd)$ satisfies the following equation:
$$\left\{\2n\ba{ll}
\ns\ds\dot P+P(A+B\Th)+(A+B\Th)^\top P+(C+D\Th)^\top P(C+D\Th)\\
\ns\ds\q+\,\Th^\top R\Th+S^\top\Th+\Th^\top S+Q=0,\qq\ae~s\in[0,T],\\
\ns\ds P(T)=G,\ea\right.$$
which is equivalent to (\ref{Riccati}).

\ms

(ii) $\Ra$ (i). Let $P(\cd)$ be the strongly regular solution of (\ref{Riccati}).
Then there exists a $\l>0$ such that
\bel{iitoi}R(s)+D(s)^\top P(s)D(s)\ges \l  I,\qq\ae~s\in[0,T].\ee
Set
$$\Th\deq-(R+D^\top PD)^{-1}(B^\top P+D^\top PC+S)\in L^2(0,T;\dbR^{m\times n}).$$
For any $u(\cd)\in\cU[0,T]$, let $X^{(u)}(\cd)$ be the solution of
$$\left\{\2n\ba{ll}
\ns\ds dX^{(u)}(s)=\big[A(s)X^{(u)}(s)+B(s)u(s)\big]ds+\big[C(s)X^{(u)}(s)+D(s)u(s)\big]dW(s),\qq s\in[0,T], \\
\ns\ds X(0)=0.\ea\right.$$
Applying It\^o's formula to $s\mapsto\langle P(s)X^{(u)}(s),X^{(u)}(s)\rangle$, we have
$$\ba{ll}
\ns\ds J^0(0,0;u(\cd))=\dbE\left\{\lan GX^{(u)}(T),X^{(u)}(T)\ran+\int_0^T\llan\begin{pmatrix}Q&S^\top\\S&R\end{pmatrix}
                                    \begin{pmatrix}X^{(u)}\\ u\end{pmatrix},
                                    \begin{pmatrix}X^{(u)}\\ u\end{pmatrix}\rran ds\right\}\\
\ns\ds=\dbE\int_0^T\[\lan\dot PX^{(u)},X^{(u)}\ran+\lan P\big(AX^{(u)}+Bu\big),X^{(u)}\ran+\lan PX^{(u)},AX^{(u)}+Bu\ran\\
\ns\ds\qq\qq~\2n+\lan P\big(CX^{(u)}+Du\big),CX^{(u)}+Du\ran+\lan QX^{(u)},X^{(u)}\ran+2\lan SX^{(u)},u\ran+\lan Ru,u\ran\]ds\\
\ns\ds=\dbE\int_0^T\[\lan\big(\dot P+PA+A^\top P+C^\top PC+Q\big)X^{(u)},X^{(u)}\ran+2\lan\big(B^\top P+D^\top PC+S\big)X^{(u)},u\ran\\
\ns\ds\qq\qq~\2n+\lan\big(R+D^\top PD\big)u,u\ran\]ds\\
\ns\ds=\dbE\int_0^T\[\lan\Th^\top\big(R+D^\top PD\big)\Th X^{(u)},X^{(u)}\ran
-2\lan\big(R+D^\top PD\big)\Th X^{(u)},u\ran+\lan\big(R+D^\top PD\big)u,u\ran\]ds\\
\ns\ds=\dbE\int_0^T\lan\big(R+D^\top PD\big)\big(u-\Th X^{(u)}\big),u-\Th X^{(u)}\ran ds. \ea$$
Noting (\ref{iitoi}) and making use of Lemma 2.5, we obtain that
$$J^0(0,0;u(\cd))=\dbE\int_0^T\lan\big(R+D^\top PD\big)\big(u-\Th X^{(u)}\big),u-\Th X^{(u)}\ran ds\ges\l\g\dbE\int_0^T|u(s)|^2ds,
\q\forall u(\cd)\in\cU[0,T],$$
for some $\g>0$. Then (i) holds. \endpf

\ms

\bf Remark 4.7. \rm From the first part of the proof of Theorem 4.6, we see that if (\ref{J>l*}) holds,
then the strongly regular solution of (\ref{Riccati}) satisfies (\ref{strong-regular})
with the same constant $\l>0$.

\ms

Combining Theorem 4.4 and Theorem 4.6, we obtain the following
corollary.

\ms

\bf Corollary 4.8. \sl Let {\rm(H1)--(H2)} and {\rm(\ref{J>l*})}
hold. Then Problem {\rm(SLQ)} is uniquely open-loop solvable at any
$(t,x)\in[0,T)\times\dbR^n$ with the open-loop optimal control
$u^*(\cd)$ being of a state feedback form:
\bel{opti-biaoshi} u^*(\cd)=-(R+D^\top PD)^{-1}(B^\top P+D^\top PC+S)X^*
-(R+D^\top PD)^{-1}(B^\top\eta+D^\top\z+D^\top P\si+\rho),\ee
where $P(\cd)$ is the unique strongly regular solution of {\rm(\ref{Riccati})}
with $(\eta(\cd),\z(\cd))$ being the adapted
solution of {\rm(\ref{eta-zeta})} and $X^*(\cd)$ being the solution
of the following closed-loop system:
\bel{closed-loop-state}\left\{\2n\ba{ll}
\ns\ds dX^*(s)=\Big\{\big[A-B(R+D^\top PD)^{-1}(B^\top P+D^\top PC+S)\big]X^*\\
\ns\ds\qq\qq\q~-B(R+D^\top PD)^{-1}(B^\top\eta+D^\top\z+D^\top P\si+\rho)+b\Big\}ds\\
\ns\ds\qq\qq~~\1n+\Big\{\big[C-D(R+D^\top PD)^{-1}(B^\top P+D^\top PC+S)\big]X^*\\
\ns\ds\qq\qq\qq~-D(R+D^\top PD)^{-1}(B^\top\eta+D^\top\z+D^\top P\si+\rho)+\si\Big\}dW(s),\qq s\in[t,T], \\
\ns\ds X^*(t)=x.\ea\right.\ee
%

\ms

\it Proof. \rm By Theorem 4.6, the Riccati equation (\ref{Riccati})
admits a unique strongly regular solution $P(\cd)\in C([0,T];\dbS^n)$.
Hence, the adapted solution $(\eta(\cd),\z(\cd))$ of {\rm(\ref{eta-zeta})}
satisfies (\ref{eta-zeta-regularity}) automatically. Now applying Theorem 4.4
and noting the remark right after Definition 2.1, we get the desired result.
\endpf

\ms

\bf Remark 4.9. \rm Under the assumptions of Corollary 4.8, when $b(\cd), \si(\cd), g, q(\cd), \rho(\cd)=0$,
the adapted solution of (\ref{eta-zeta}) is $(\eta(\cd),\z(\cd))\equiv(0,0)$. Thus, for Problem ${\rm(SLQ)}^0$,
the unique optimal control $u^*(\cd)$ at initial pair $(t,x)\in[0,T)\times\dbR^n$ is given by
\bel{opti-biaoshi-0} u^*(\cd)=-(R+D^\top PD)^{-1}(B^\top P+D^\top PC+S)X^*,\ee
with $P(\cd)$ being the unique strongly regular solution of {\rm(\ref{Riccati})}
and $X^*(\cd)$ being the solution of the following closed-loop system:
\bel{closed-loop-state-0}\left\{\2n\ba{ll}
\ns\ds dX^*(s)=\big[A-B(R+D^\top PD)^{-1}(B^\top P+D^\top PC+S)\big]X^*ds\\
\ns\ds\qq\qq~~\1n+\big[C-D(R+D^\top PD)^{-1}(B^\top P+D^\top PC+S)\big]X^*dW(s),\qq s\in[t,T], \\
\ns\ds X^*(t)=x.\ea\right.\ee
Moreover, by (\ref{Value}), the value function of Problem ${\rm(SLQ)}^0$ is given by
\bel{} V^0(t,x)=\langle P(t)x,x\rangle, \qq (t,x)\in[0,T]\times\dbR^n.\ee

\section{Finiteness of Problem (SLQ) and Convexity of Cost Functional}

We have seen that the uniform convexity of the cost functional
implies the open-loop and closed-loop solvabilities of Problem
(SLQ). We expect that the finiteness of Problem (SLQ) should be
closely related to the convexity of the cost functional. A main
purpose of this section is to make this clear. Other relevant issues
will also be discussed. First, we introduce the following:
\bel{L}\L(s,P(\cd))=\begin{pmatrix}\scriptstyle\dot
P(s)+P(s)A(s)+A(s)^\top P(s)+C(s)^\top P(s)C(s)+Q(s)
&\scriptstyle P(s)B(s)+C(s)^\top P(s)D(s)+S(s)^\top\\
\scriptstyle B(s)^\top P(s)+D(s)^\top P(s)C(s)+S(s)&\scriptstyle
R(s)+D(s)^\top P(s)D(s)\end{pmatrix},\ee
for any $P(\cd)\in AC(t,T;\dbS^n)$ which is the set of all absolutely continuous
functions $P:[t,T]\to\dbS^n$. Let
$$\cP[t,T]=\Big\{P(\cd)\in
AC(t,T;\dbS^n)\bigm|P(T)\les G,~\L(s,P(\cd))\ges0,~\ae
s\in[t,T]\Big\}.$$
We have the following result.

\ms

\bf Proposition 5.1. \sl Let {\rm(H1)--(H2)} hold, and $t\in[0,T)$
be given. Among the following statements:

\ms

{\rm(i)} Problem {\rm(SLQ)} is finite at $t$.

\ms

{\rm(ii)} Problem {$\rm(SLQ)^0$} is finite at $t$.

\ms

{\rm(iii)} There exists a $P(t)\in\dbS^n$ such that
\bel{V0}V^0(t,x)=\langle P(t)x,x\rangle,\qq\forall x\in\dbR^n.\ee

\ms

{\rm(iv)} The map $u(\cd)\mapsto J(t,x;u(\cd))$ is convex, for any
$x\in\dbR^n$.

\ms

{\rm(v)} $\cP[t,T]\ne\emptyset$.

\ms

\no the following implications hold:

\ms

\rm

(i) $\Ra$ (ii) $\Ra$ (iii) $\Ra$ (iv); (v) $\Ra$ (ii).

\ms

\it Proof. \rm (i) $\Ra$ (ii). By Proposition 3.1, for any
$x\in\dbR^n$ and $u(\cd)\in\cU[t,T]$, we have
$$\ba{ll}
\ns\ds V(t,x)+V(t,-x)\les J(t,x;u(\cd))+J(t,-x;-u(\cd))\\
\ns\ds\qq\qq\qq\qq\,=2\[\lan M_2(t)u,u\ran+2\lan M_1(t)x,u\ran+\lan M_0(t)x,x\ran+c_t\]\\
\ns\ds\qq\qq\qq\qq\,=2\[J^0(t,x;u(\cd))+c_t\],\ea$$
which implies (ii).

\ms

(ii) $\Ra$ (iii) can be shown by a simple adoption of the well-known
result in the deterministic case (see \cite{Clements-Anderson-Moylan
1977, Anderson-Moore 1989}).

\ms

(iii) $\Ra$ (iv). By Corollary 3.3, if $u(\cd)\mapsto J(t,x;u(\cd))$
is not convex, then $J^0(t,0;u(\cd))<0$ for some
$u(\cd)\in\cU[t,T]$. By Corollary 3.2, we have
$$J^0(t,x;\l u(\cd))=J^0(t,x;0)+\l^2J^0(t,0;u(\cd))+\l\dbE\int_t^T\lan\cD
J^0(t,x;0)(s),u(s)\ran ds,\qq\forall \l\in\dbR^n.$$
Letting $\l\to\i$, we obtain
$$V^0(t,x)\les\lim_{\l\to\i}J^0(t,x;\l u(\cd))=-\i,$$
which is a contradiction.

\ms

(v) $\Ra$ (ii). For any $(t,x)\in[0,T)\times\dbR^n$,
$u(\cd)\in\cU[t,T]$, and any $P(\cd)\in AC(t,T;\dbS^n)$, one
has
$$\ba{ll}
\ns\ds\dbE\lan P(T)X(T),X(T)\ran-\lan
P(t)x,x\ran\\
\ns\ds=\dbE\int_t^T\bigg\{\llan\[\dot P(s)+P(s)A(s)+A(s)^\top
P(s)+C(s)^\top P(s)C(s)\]X(s),X(s)\rran\\
\ns\ds\qq\qq~+2\llan\[B(s)^\top P(s)+D(s)^\top
P(s)C(s)\]X(s),u(s)\rran+\lan D(s)^\top
P(s)D(s)u(s),u(s)\ran\bigg\}ds.\ea$$
Hence, if $\cP[t,T]\ne\emptyset$, then by taking
$P(\cd)\in\cP[t,T]$, one has
$$\ba{ll}
\ns\ds J^0(t,x;u(\cd))=\dbE\left\{\lan
GX(T),X(T)\ran+\int_t^T\llan\begin{pmatrix}Q(s)&S(s)^\top\\
S(s)&R(s)\end{pmatrix}\begin{pmatrix}X(s)\\
u(s)\end{pmatrix},\begin{pmatrix}X(s)\\ u(s)\end{pmatrix}\rran
ds\right\}\\
\ns\ds=\lan
P(t)x,x\ran+\dbE\left\{\lan\big[G-P(T)\big]X(T),X(T)\ran+\int_t^T\llan\L(s,P(\cd))\begin{pmatrix}X(s)\\
u(s)\end{pmatrix},\begin{pmatrix}X(s)\\
u(s)\end{pmatrix}\rran ds\right\}\\
\ns\ds\ges\lan P(t)x,x\ran,\qq\forall u(\cd)\in\cU[t,T].\ea$$
This implies that the corresponding Problem (SLQ)$^0$ is finite at $t$.
\endpf

\ms

It is worth to point out that the convexity of the map
$u(\cd)\mapsto J^0(t,x;u(\cd))$ is not sufficient for the finiteness
of Problem $\hb{(SLQ)}^0$. We present the following example (see
also \cite{Mou-Yong 2006} for an example of a quadratic functional
in Hilbert space).

\ms

\bf Example 5.2. \rm Consider the following one-dimensional
controlled SDE:
\bel{countEx-1}\left\{\2n\ba{ll}
\ns\ds dX(s)=u(s)ds+X(s)dW(s),\q s\in[t,1],\\
\ns\ds X(t)=x,\ea\right.\ee
and the cost functional:
\bel{countEx-2}J^0(t,x;u(\cd))=\dbE\left[-X(1)^2+\int_t^1
e^{1-s}u(s)^2ds\right].\ee
We claim that
\bel{countEx-3}J^0(0,0;u(\cd))\ges0,\qq \forall
u(\cd)\in\cU[0,T],\ee
which, by Corollary 3.3, is equivalent to the convexity of
$u(\cd)\mapsto J^0(0,x;u(\cd))$, but
\bel{countEx-4}V^0(0,x)=-\i,\qq \forall x\neq 0.\ee
To show the above, let $u(\cd)\in\cU[0,T]$ and $X(\cd)\equiv
X(\cd\,;0,x,u(\cd))$ be the solution of (\ref{countEx-1}) with
$t=0$. By the variation of constants formula,
$$X(s)=xe^{W(s)-{1\over2}s}+e^{W(s)-{1\over2}s}\int_0^se^{{1\over2}r-W(r)}u(r)dr,
\qq s\in[0,1].$$
Taking $x=0$ and noting that $e^{2[W(1)-W(r)]-(1-r)}$ is independent
of $\cF_r$, we have
$$\ba{ll}
\ns\ds\dbE\big[X(1)^2\big]=\dbE\left[\int_0^1e^{W(1)-W(r)-{1\over2}(1-r)}u(r)dr\right]^2
\les\dbE\int_0^1e^{2[W(1)-W(r)]-(1-r)}u(r)^2dr\\
\ns\ds\qq\qq\1n~=\int_0^1\dbE e^{2[W(1)-W(r)]-(1-r)}\dbE[u(r)^2]dr
=\int_0^1e^{1-r}\dbE[u(r)^2]dr=\dbE\int_0^1e^{1-r}u(r)^2dr,\ea$$
and hence
$$J^0(0,0;u(\cd))=\dbE\left[-X(1)^2+\int_0^1 e^{1-s}u(s)^2ds\right]\ges0,\qq\forall u(\cd)\in\cU[0,T].$$
On the other hand, taking $x\ne0$ and $u(s)=\l
e^{W(s)-{1\over2}s},~\l\in\dbR$, we have
$$X(1)=(x+\l)e^{W(1)-{1\over2}}.$$
Therefore,
$$\ba{ll}
\ns\ds J^0(0,x;u(\cd))=\dbE\left[-X(1)^2+\int_0^1
e^{1-s}u(s)^2ds\right]
=-\dbE\[(x+\l)^2e^{2W(1)-1}\]+\l^2\dbE\int_0^1 e^{1-s}e^{2W(s)-s}ds\\
\ns\ds\qq\qq\q~~\1n=-(x+\l)^2e+\l^2e=-(x^2+2\l x)e.\ea$$
Letting $|\l|\to\i$, with $\l x>0$, in the above, we obtain
$V^0(0,x)=-\i$. This proves our claim.

\ms

The above example tells us that, besides the convexity of
$u(\cd)\mapsto J^0(t,x;u(\cd))$, one needs some additional
condition(s) in order to get the finiteness of Problem (SLQ)$^0$ at
$t$. To find such a condition, let us make some observations.
Suppose $u(\cd)\mapsto J^0(0,x;u(\cd))$ is convex, which, by
Corollary 3.3, is equivalent to the following:
\bel{t=0-convex}J^0(0,0;u(\cd))\ges0,\qq\forall u(\cd)\in\cU[0,T].\ee
Then for any $\e>0$, consider state equation (\ref{state}) (with
$b(\cd), \si(\cd)=0$) and the following cost functional:
\bel{cost-e}\ba{ll}
\ns\ds J^0_\e(t,x;u(\cd))\deq\dbE\left\{\langle GX(T),X(T)\rangle
+\int_t^T\llan\begin{pmatrix}Q(s)&S(s)^\top\\
                                          S(s)&R(s)+\e I\end{pmatrix}\begin{pmatrix}
                                          X(s)\\ u(s)\end{pmatrix},
                                          \begin{pmatrix}
                                          X(s)\\
                                          u(s)\end{pmatrix}\rran ds\right\}\\
\ns\ds\qq\qq\qq\2n=J^0(t,x;u(\cd))+\e\dbE\int_0^T|u(s)|^2ds.\ea\ee
Denote the corresponding optimal control problem and value function
by Problem ${\rm(SLQ)}^0_\e$ and $V^0_\e(\cd\,,\cd)$, respectively.
By Corollary 3.3 and the convexity of $u(\cd)\mapsto
J^0(0,0;u(\cd))$, one has
$$J^0_\e(0,0;u(\cd))=J^0(0,0;u(\cd))+\e\dbE\int_0^T|u(s)|^2ds\ges\e\dbE\int_0^T|u(s)|^2ds,\qq\forall u(\cd)\in\cU[0,T],$$
i.e., $u(\cd)\mapsto J^0_\e(0,0;u(\cd))$ is uniformly convex. Hence,
it follows from Theorem 4.6 that the following Riccati equation:
\bel{Ric-e}\left\{\2n\ba{ll}
\ns\ds\dot P_\e+P_\e A+A^\top P_\e+C^\top P_\e C+Q\\
\ns\ds\q\2n~-(P_\e B+C^\top\1n P_\e D+S^\top)(R+\e I+D^\top\1n P_\e D)^{-1}(B^\top P_\e+D^\top\1n P_\e C+S)=0,\qq \ae~s\in[0,T],\\
\ns\ds P_\e(T)=G,\ea\right.\ee
admits a unique strongly regular solution $P_\e(\cd)\in
C([0,T];\dbS^n)$ such that (noting Remark 4.7)
\bel{R+DPD>e}R(t)+\e I+D(t)^\top P_\e(t)D(t)\ges\e I,\qq\ae
t\in[0,T].\ee
Now, we are ready to state and prove the following result which is a
characterization of the finiteness of Problem (SLQ)$^0$.

\ms

\bf Theorem 5.3. \sl Let {\rm(H1)--(H2)} and {\rm(\ref{t=0-convex})}
hold. For any $\e>0$, let $P_\e(\cd)$ be the unique strongly regular
solution of the Riccati equation {\rm(\ref{Ric-e})}. Then Problem
${\rm(SLQ)}^0$ is finite if and only if $\{P_\e(0)\}_{\e>0}$ is
bounded from below. In this case, the limit
\bel{lim P_e}\lim_{\e\to0}P_\e(t)=P(t),\qq\forall t\in[0,T],\ee
exists, and {\rm(\ref{V0})} holds. Moreover,
\bel{R+DPD>0}R(t)+D(t)^\top P(t)D(t)\ges0,\qq\ae~t\in[0,T],\ee
and
\bel{N<P<M}N(t)\les P(t)\les M_0(t),\qq\forall t\in[0,T],\ee
where $M_0(\cd)$ is the solution to the Lyapunov equation
{\rm(\ref{3.8}), and
\bel{N}N(t)=\big[\F_A(t)^\top\big]^{-1}\left\{P(0)-\int_0^t\F_A(s)^\top
\[C(s)^\top M_0(s)C(s)+Q(s)\]\F_A(s)ds\right\}\F_A(t)^{-1},\ee
with $\F_A(\cd)$ being the solution to the following:
\bel{FA5.5}\left\{\2n\ba{ll}
\ns\ds\dot\F_A(s)=A(s)\F_A(s),\qq s\ges0,\\
\ns\ds\F_A(0)=I.\ea\right.\ee
In particular, if Problem ${\rm(SLQ)}^0$ is finite at $t=0$, then it is finite.

\ms

\it Proof. {\it Necessity}. \rm Suppose Problem $\hb{(SLQ)}^0$ is finite and let
$P:[0,T]\to\dbS^n$ such that {\rm(\ref{V0})} holds. For any $\e_2>\e_1>0$, we have
$$J^0_{\e_2}(t,x;u(\cd))\ges J^0_{\e_1}(t,x;u(\cd))\ges J^0(t,x;u(\cd)),\qq\forall(t,x)\in[0,T]\times\dbR^n,\q\forall u(\cd)\in\cU[t,T].$$
Hence (noting Remark 4.9),
$$\ba{ll}
\ns\ds\langle P_{\e_2}(t)x,x\rangle=V^0_{\e_2}(t,x)=\inf_{u(\cd)\in\cU[t,T]}J^0_{\e_2}(t,x;u(\cd))\\
\ns\ds\qq\qq\q~\2n\ges\inf_{u(\cd)\in\cU[t,T]}J^0_{\e_1}(t,x;u(\cd))=V^0_{\e_1}(t,x)=\langle P_{\e_1}(t)x,x\rangle\\
\ns\ds\qq\qq\q~\2n\ges\inf_{u(\cd)\in\cU[t,T]}J^0(t,x;u(\cd))=V^0(t,x)=\langle P(t)x,x\rangle,
\qq\forall(t,x)\in[0,T]\times\dbR^n.\ea$$
Thus $\{P_\e(t)\}_{\e>0}$ is a nondecreasing sequence with lower bound $P(t)$ and therefore
has a limit $\bar P(t)$ with
\bel{barP>P}\bar P(t)\equiv\lim_{\e\to0}P_\e(t)\ges P(t),\qq\forall t\in[0,T].\ee
On the other hand, for any $\d>0$, we can find a $u^\d(\cd)\in\cU[t,T]$, such that
$$\ba{ll}
\ns\ds V^0_\e(t,x)\les
J^0(t,x;u^\d(\cd))+\e\dbE\int_t^T|u^\d(s)|^2ds\les
V^0(t,x)+\d+\e\dbE\int_t^T|u^\d(s)|^2ds.\ea$$
Letting $\e\to 0$, we obtain that
$$\langle \bar P(t)x,x\rangle=\lim_{\e\to 0}\langle P_\e(t)x,x\rangle
=\lim_{\e\to 0}V^0_\e(t,x)\les V^0(t,x)+\d=\langle P(t)x,x\rangle+\d,$$
from which we see that
\bel{barP<P}\langle \bar P(t)x,x\rangle\les\langle P(t)x,x\rangle
\qq\forall (t,x)\in[0,T]\times\dbR^n.\ee
Combining (\ref{barP>P})--(\ref{barP<P}), we obtain (\ref{lim P_e}) with $P(\cd)$ satisfying (\ref{V0}).
Moreover, letting $\e\to0$ in (\ref{R+DPD>e}), we obtain (\ref{R+DPD>0}).

\ms

{\it Sufficiency}. Suppose there exists a $\b\in\dbR$ such
that
$$P_\e(0)\ges\b I,\qq\forall\e>0,$$
then for any $x\in\dbR^n$ and $u(\cd)\in\cU[0,T]$, we have
$$J^0(0,x;u(\cd))+\e\dbE\int_0^T|u(s)|^2ds\ges V^0_\e(0,x)=\langle P_\e(0)x,x\rangle\ges \b|x|^2,\qq\forall\e>0.$$
Letting $\e\to0$ in the above, we obtain
$$J^0(0,x;u(\cd))\ges \b|x|^2,\qq\forall x\in\dbR^n,\q\forall u(\cd)\in\cU[0,T],$$
which implies the finiteness of Problem $\hb{(SLQ)}^0$ at $t=0$.

\ms

Now, let $P(0)\in\dbS^n$ such that $V^0(0,x)=\langle P(0)x,x\rangle,\forall x\in\dbR^n$.
Then
\bel{p0<pe0}P(0)\les P_\e(0),\qq\forall \e>0.\ee
Also, by Remark 4.9 and Proposition 3.1,
$$\langle P_\e(t)x,x\rangle=V^0_\e(t,x)\les
J^0_\e(t,x;0)=J^0(t,x;0)=\langle M_0(t)x,x\rangle,\qq\forall
x\in\dbR^n.$$
This leads to
\bel{5.16}P_\e(t)\les M_0(t),\qq t\in[0,T],\qq\forall\e>0.\ee
On the other hand, let $\F_A(\cd)$ be the solution of (\ref{FA5.5}) and set
$$\Pi_\e\deq(P_\e B+C^\top P_\e D+S^\top)(R_\e+D^\top P_\e D)^{-1}(B^\top P_\e+D^\top P_\e C+S)\ges0.$$
Differentiating $s\mapsto\F_A(s)^\top P_\e(s)\F_A(s)$, we obtain
$${d\over ds}\[\F_A(s)^\top P_\e(s)\F_A(s)\]=\F_A(s)^\top\[\Pi_\e(s)-C(s)^\top
P_\e(s)C(s)-Q(s)\]\F_A(s).$$
Hence, combining (\ref{p0<pe0})--(\ref{5.16}), we have
$$\ba{ll}
\ns\ds\F_A(t)^\top
P_\e(t)\F_A(t)=P_\e(0)+\int_0^t\F_A(s)^\top\[\Pi_\e(s)-C(s)^\top P_\e(s)C(s)-Q(s)\]\F_A(s)ds\\
\ns\ds\qq\qq\qq\qq~\ges P(0)-\int_0^t\F_A(s)^\top\[C(s)^\top
M_0(s)C(s)+Q(s)\]\F_A(s)ds.\ea$$
Thus,
\bel{5.17}\ba{ll}
\ns\ds P_\e(t)\ges\big[\F_A(t)^\top\big]^{-1}\left\{P(0)-\int_0^t\F_A(s)^\top \[C(s)^\top M_0(s)C(s)+Q(s)\]\F_A(s)ds\right\}\F_A(t)^{-1}\\
\ns\ds\qq\1n~\equiv N(t),\qq t\in[0,T].\ea\ee
Then using the same argument as in the previous paragraph, we can show that
\bel{youxian-0}J^0(t,x;u(\cd))\ges\langle N(t)x,x\rangle,\qq\forall (t,x)\in[0,T]\times\dbR^n,\q\forall u(\cd)\in\cU[t,T],\ee
which implies the finiteness of Problem $\hb{(SLQ)}^0$. Moreover, let
$P:[0,T]\to\dbS^n$ such that {\rm(\ref{V0})} holds, then
$$\langle N(t)x,x\rangle\les\inf_{u(\cd)\in\cU[t,T]}J^0(t,x;u(\cd))
=\langle P(t)x,x\rangle\les J^0(t,x;0)=\langle M_0(t)x,x\rangle,\qq\forall (t,x)\in[0,T]\times\dbR^n,$$
and (\ref{N<P<M}) follows.

\ms

Finally, if Problem ${\rm(SLQ)}^0$ is finite at $t=0$, then (\ref{p0<pe0}) holds
and the finiteness of Problem ${\rm(SLQ)}^0$ therefore follows.
\endpf

\ms

The following is another sufficient condition for the finiteness of
Problem (SLQ)$^0$, which is a corollary of Theorem 5.3 and Proposition 5.1, (v)
$\Ra$ (ii).

\ms

\bf Corollary 5.4. \sl Let {\rm(H1)--(H2)} hold. If there exists a
$\D(\cd)\in L^1(0,T;\dbS^n_+)$ such that
\bel{Th3.6-2} R+D^\top PD\ges \big(B^\top P+D^\top P
C+S\big)\D^{-1}\big(PB+C^\top PD+S^\top\big),\qq\ae~s\in[0,T],\ee
where $P(\cd)$ is the solution of the following Lyapunov equation:
\bel{Th3.6-1}\left\{\2n\ba{ll}
\ns\ds \dot P(s)+P(s)A(s)+A(s)^\top P(s)+C(s)^\top P(s)C(s)+Q(s)=\D(s),\qq \ae~s\in [0,T],\\
\ns\ds P(T)\les G.\ea\right.\ee
Then Problem {\rm(SLQ)$^0$} is finite.

\ms

\it Proof. \rm Under our condition, by Lemma 2.6, one has
$$\L(s,P(\cd))=\begin{pmatrix}\scriptstyle\D(s)
&\scriptstyle P(s)B(s)+C(s)^\top P(s)D(s)+S(s)^\top\\
\scriptstyle B(s)^\top P(s)+D(s)^\top P(s)C(s)+S(s)&\scriptstyle
R(s)+D(s)^\top P(s)D(s)\end{pmatrix}\ges0,\qq\ae~s\in[0,T].$$
Hence, $P(\cd)\in\cP[0,T]$ and Problem (SLQ)$^0$ is finite at $t=0$.
Then the finiteness of Problem (SLQ)$^0$ follows from Theorem 5.3 immediately.
\endpf

\ms

We now return to the study of convexity of the map $u(\cd)\mapsto
J^0(t,0;u(\cd))$. First, from the definition of $M_2(t)$ (see
(\ref{J-rep})), we see that $M_2(t)\ges0$ if and only if
\bel{3.20}R(\cd)\ges-\(\h L_t^*G\h
L_t+L^*_tQL_t+SL_t+L^*_tS^\top\),\ee
with the right hand side being non-positive. Thus, unlike the
well-known situation for the deterministic LQ problems (for which
$R(\cd)\ges0$ is necessary for $M_2(t)\ges0$ (\cite{Yong-Zhou
1999})), $R(\cd)$ does not have to be positive semi-definite.
Actually, as shown by examples in \cite{Chen-Li-Zhou 1998, Yong-Zhou
1999}, $R(\cd)$ could even be negative definite in some extent. Let
us now take a closer look at this issue.

\ms

Note that when $u(\cd)\mapsto J^0(0,0;u(\cd))$ is convex, for any $\e>0$, the
unique strongly regular solution $P_\e(\cd)$ to the Riccati equation (\ref{Ric-e})
satisfies (\ref{R+DPD>e}) and (\ref{5.16}). Hence,
\bel{R+DMD>0}R(t)+D(t)^\top M_0(t)D(t)\ges0,\qq\ae t\in[0,T],\ee
or equivalently,
\bel{R+D[M]D>0}\ba{ll}
\ns\ds R(t)+D(t)^\top\dbE\bigg\{\big[\F(T)\F(t)^{-1}\big]^\top
G\big[\F(T)\F(t)^{-1}\big]\\
\ns\ds\qq\qq\qq\qq+\int_t^T\big[\F(s)\F(t)^{-1}\big]^\top
Q(s)\big[\F(s)\F(t)^{-1}\big]ds\bigg\}D(t)\ges0,\qq\ae
t\in[0,T].\ea\ee
This is another necessary condition for the finiteness of Problem
(SLQ)$^0$, which is easier to check. From (\ref{R+D[M]D>0}), we see
that if $R(\cd)$ happens to be negative definite, then in order
$u(\cd)\mapsto J^0(0,0;u(\cd))$ to be convex, it is necessary that
$D(\cd)$ is injective, and either $G$ or $Q(\cd)$ (or both) has to
be positive enough to compensate. Note that $D(s)$ was assumed to be
invertible in \cite{Qian-Zhou 2013}. Therefore, in some sense, our
result justifies the assumption of \cite{Qian-Zhou 2013}.

\ms

The following gives a little improvement when more restrictive
conditions are assumed.

\ms

\bf Proposition 5.5. \sl Let {\rm(H1)--(H2)} hold. Suppose that
\bel{B=0}B(\cd)=0,\qq C(\cd)=0,\qq S(\cd)=0.\ee
Then the map $u(\cd)\mapsto J^0(0,0;u(\cd))$ is convex if and only
if {\rm(\ref{R+D[M]D>0})} holds. In this case, Problem
{\rm(SLQ)$^0$} is closed-loop solvable.

\ms

\it Proof. \rm It suffices to prove the sufficiency. Note that in
the current case, the corresponding Riccati equation becomes
\bel{}\left\{\2n\ba{ll}
\ns\ds\dot P(s)+P(s)A(s)+A(s)^\top P(s)+Q(s)=0,\qq \ae~s\in[0,T],\\
\ns\ds P(T)=G,\ea\right.\ee
whose solution is $M_0(\cd)$. If {\rm(\ref{R+D[M]D>0})} holds, then it is easy to verify
that $M_0(\cd)$ is regular. Consequently, by Theorem 4.4, Problem (SLQ)$^0$ is closed-loop
solvable, and hence $u(\cd)\mapsto J^0(0,0;u(\cd))$ is convex. \endpf

\ms

Note that in the case (\ref{B=0}), we have
$$M_0(t)=\big[\F_A(T)\F_A(t)^{-1}\big]^\top
G\big[\F_A(T)\F_A(t)^{-1}\big]+\int_t^T\big[\F_A(s)\F_A(t)^{-1}\big]^\top
Q(s)\big[\F_A(s)\F_A(t)^{-1}\big]ds,$$
with $\F_A(\cd)$ being the solution of (\ref{FA5.5}). Hence,
(\ref{R+D[M]D>0}) can also be written as
\bel{5.7}\ba{ll}
\ns\ds R(t)+D(t)^\top\bigg\{\big[\F_A(T)\F_A(t)^{-1}\big]^\top
G\big[\F_A(T)\F_A(t)^{-1}\big]\\
\ns\ds\qq\qq\qq\q+\int_t^T\big[\F_A(s)\F_A(t)^{-1}\big]^\top
Q(s)\big[\F_A(s)\F_A(t)^{-1}\big]ds\bigg\}D(t)\ges0,\qq\ae~t\in[0,T].\ea\ee

\ms

As we pointed out earlier, Problem $\hb{(SLQ)}^0$ may still be
infinite at some initial pair $(t,x)$ even if the cost functional is
convex. In this case, by Theorem 5.3, the sequence
$\{P_\e(t)\}_{\e>0}$ diverges. The following result is concerned
with the divergence speed of  $\{P_\e(t)\}_{\e>0}$.

\ms

\bf Proposition 5.6. \sl Let {\rm(H1)--(H2)} and
{\rm(\ref{t=0-convex})} hold. For any $\e>0$, let $P_\e(\cd)$ be the
unique strongly regular solution of {\rm(\ref{Ric-e})}. Then for any
$\a\ges 1$ and $t\in[0,T]$, the sequence $\{\e^\a P_\e(t)\}_{\e>0}$
converges.

\ms

\it Proof. \rm Let $\e>0$ and consider Problem $\hb{(SLQ)}_\e^0$. By Proposition 3.1 (noting $M_2(t)\ges0$),
$$\ba{ll}
\ns\ds J_\e^0(t,x;u(\cd))=J^0(t,x;u(\cd))+\e\dbE\int_t^T|u(s)|^2ds\\
\ns\ds\qq\qq\q~=\lan [M_2(t)+\e I] u,u\ran+2\lan M_1(t)x,u\ran+\lan M_0(t)x,x\ran\\
\ns\ds\qq\qq\q~=\Big|[M_2(t)+\e I]^{1\over2}u+[M_2(t)+\e
I]^{-{1\over2}}M_1(t)x\Big|^2\\
\ns\ds\qq\qq\qq~~+\lan M_0(t)x,x\ran-\lan[M_2(t)+\e
I]^{-1}M_1(t)x,M_1(t)x\ran.\ea$$
Thus
\bel{}\langle P_\e(t)x,x\rangle=\inf_{u(\cd)\in\cU[t,T]}J_\e^0(t,x;u(\cd))
=\lan M_0(t)x,x\ran-\lan[M_2(t)+\e I]^{-1}M_1(t)x,M_1(t)x\ran.\ee
For any $\a\ges1$, we have
$$0\les\e^\a\lan[M_2(t)+\e I]^{-1}M_1(t)x,M_1(t)x\ran\les \e^{\a-1}\lan M_1(t)x,M_1(t)x\ran
\les\lan M_1(t)x,M_1(t)x\ran,\qq\forall 0<\e\les1.$$
Consequently, for any $x\in\dbR^n$, $\{\e^\a \langle P_\e(t)x,x\rangle\}_{\e>0}$
has a convergent subsequence, and hence the sequence $\{\e^\a \langle P_\e(t)x,x\rangle\}_{\e>0}$
itself converges as $\e\to0$ since it is nondecreasing. The result therefore follows.
\endpf

\ms

From Remark 4.7 and 4.9, we see that if the uniformly convex condition
(\ref{J>l*}) holds, then Problem $\hb{(SLQ)}^0$ is finite and
\bel{uni-posi}R(s)+D(s)^\top P(s)D(s)\ges \l I, \qq\ae~s\in[0,T],\ee
where $P:[0,T]\to\dbS^n$ is the function such that {\rm(\ref{V0})} holds. The following
result shows that the converse is also true.

\ms

\bf Theorem 5.7. \sl Let {\rm(H1)--(H2)} hold. Suppose Problem
${\rm(SLQ)}^0$ is finite and let $P:[0,T]\to\dbS^n$ such that
{\rm(\ref{V0})} holds. If {\rm(\ref{uni-posi})} holds for some
$\l>0$, then $P(\cd)$ solves the Riccati equation
{\rm(\ref{Riccati})}. Consequently, the map $u(\cd)\mapsto
J^0(0,0;u(\cd))$ is uniformly convex.

\ms

\it Proof. \rm  For any $\e>0$, let $P_\e(\cd)$ be the unique
strongly regular solution of {\rm(\ref{Ric-e})}. By Theorem 5.3,
$$P_\e(t)\searrow P(t),\qq\hb{as}\q\e\searrow0,\q\forall t\in[0,T].$$
Note that $P_\e(\cd)\les M_0(\cd)~\forall \e>0$ and by (\ref{N<P<M}), $P(\cd)$ is bounded.
Thus, $\{P_\e(t)\}_{\e>0}$ is uniformly bounded. Also, we have
$$R(s)+D(s)^\top P_\e(s)D(s)\ges R(s)+D(s)^\top P(s)D(s)\ges \l I, \qq\ae~s\in[0,T],\q\forall \e>0.$$
Then it follows from the dominated convergence theorem that
$$P_\e A+A^\top P_\e+C^\top P_\e C+Q-(P_\e B+C^\top P_\e D+S^\top)(R+\e I+D^\top P_\e D)^{-1}(B^\top P_\e+D^\top P_\e C+S)\equiv\L_\e $$
converges to
$$PA+A^\top P+C^\top PC+Q-(PB+C^\top PD+S^\top)(R+D^\top PD)^{-1}(B^\top P+D^\top PC+S)\equiv\L$$
in $L^1$ as $\e\to 0$. Therefore,
$$P(t)=\lim_{\e\to0}P_\e(t)=G+\lim_{\e\to0}\int_t^T\L_\e(s)ds=G+\int_t^T\L(s)ds,$$
which, together with (\ref{uni-posi}), implies that $P(\cd)$ is a strongly regular solution of {\rm(\ref{Riccati})}.
Consequently, by Theorem 4.6, $u(\cd)\mapsto J^0(0,0;u(\cd))$ is uniformly convex.
\endpf

\ms

We now look at the following case:
\bel{D=0}D(\cd)=0,\qq R(\cd)\gg0.\ee
Note that although $D(\cd)=0$, since $C(\cd)$ is not necessarily
zero, our state equation is still an SDE. For such a case, the above
results can be restated as follows.

\ms

\bf Theorem 5.8. \sl Let {\rm(H1)--(H2)} and {\rm(\ref{D=0})} hold.
Then the following statements are equivalent:

\ms

{\rm(i)} Problem {\rm(SLQ)} is finite at $t=0$;

\ms

{\rm(ii)} Problem ${\rm(SLQ)}^0$ is finite at $t=0$;

\ms

{\rm(iii)} The map $u(\cd)\mapsto J^0(0,0;u(\cd))$ is uniformly convex;

\ms

{\rm(iv)} The following Riccati equation
\bel{Ric-D=0}\left\{\2n\ba{ll}
\ns\ds\dot P+PA+A^\top P+C^\top PC+Q-(PB+S^\top)R^{-1}(B^\top P+S)=0,\qq \ae~s\in[0,T],\\
\ns\ds P(T)=G,\ea\right.\ee
admits a unique solution $P(\cd)\in C([0,T];\dbS^n)$;

\ms

{\rm(v)} Problem {\rm(SLQ)} is uniquely closed-loop solvable;

\ms

{\rm(vi)} Problem {\rm(SLQ)} is uniquely open-loop solvable.

\ms

\it Proof. \rm (i) $\Ra$ (ii) follows from Proposition 5.1.

\ms

(ii) $\Ra$ (iii). By Theorem 5.3, Problem ${\rm(SLQ)}^0$ is finite. Since $D(\cd)=0, R(\cd)\gg0$,
(\ref{uni-posi}) holds for some $\l>0$, and the result follows from Theorem 5.7.

\ms

(iii) $\Leftrightarrow$ (iv). In the case of (\ref{D=0}), the corresponding Riccati equation becomes (\ref{Ric-D=0}).
If $P(\cd)\in C([0,T];\dbS^n)$ is a solution of (\ref{Ric-D=0}), then it is automatically strongly regular.
Thus, by Theorem 4.6, we obtain the equivalence of (iii) and (iv).

\ms

(iii) $\Ra$ (vi) follows from Proposition 4.1 and (iv) $\Ra$ (v) follows from Theorem 4.4.

\ms

Finally, (v) $\Ra$ (i) and (vi) $\Ra$ (i) are trivial.
\endpf

\ms

An interesting point of the above is that under condition
(\ref{D=0}), finiteness of Problem (SLQ) implies the closed-loop
solvability of Problem (SLQ). In the deterministic case, such a fact
was firstly revealed in \cite{Zhang 2005} for two-person zero-sum
differential games, and was proved in \cite{Yong 2015} for
deterministic LQ problems by means of Fredholm operators.

\section{Minimizing Sequences and Open-Loop Solvabilities}

In Section 4, we showed that under the uniform convexity condition (\ref{J>l*}),
Problem (SLQ) is open-loop solvable and the open-loop optimal control has a linear
state feedback representation. In this section, we study the open-loop solvability
of Problem (SLQ) without the the uniform convexity condition.

\ms

First we construct a minimizing sequence for Problem (SLQ) when it
is finite.

\ms

\bf Theorem 6.1. \sl Let {\rm(H1)--(H2)} hold. Suppose Problem {\rm(SLQ)} is finite. For any $\e>0$,
let $P_\e(\cd)$ be the unique strongly regular solution to the Riccati equation {\rm(\ref{Ric-e})}.
Further, let $(\eta_\e(\cd), \z_\e(\cd))$ and $X_\e(\cd)\equiv X_\e(\cd\,;t,x)$ be the (adapted) solutions
to the following BSDE and closed-loop system, respectively:
\bel{eta-zeta-e}\left\{\2n\ba{ll}
\ns\ds d\eta_\e(s)=-\[(A+B\Th_\e)^\top\eta_\e+(C+D\Th_\e)^\top\z_\e
+(C+D\Th_\e)^\top P_\e\si-\Th_\e^\top\rho+P_\e b+q\]ds\\
\ns\ds\qq\qq\qq\qq\qq\qq\qq\qq\qq\qq~+\z_\e dW(s),\qq s\in[0,T],\\
\ns\ds\eta_\e(T)=g,\ea\right.\ee
\bel{close-state-e}\left\{\2n\ba{ll}
\ns\ds dX_\e(s)=\[(A+B\Th_\e)X_\e+Bv_\e+b\]ds+\[(C+D\Th_\e)X_\e+Dv_\e+\si\] dW(s),\qq s\in[t,T], \\
\ns\ds X_\e(t)=x,\ea\right.\ee
where
\bel{Th-e}\left\{\2n\ba{ll}
\ns\ds \Th_\e=-\big(R+\e I+D^\top P_\e D\big)^{-1}\big(B^\top P_\e+D^\top P_\e C+S\big),\\
\ns\ds v_\e=-\big(R+\e I+D^\top P_\e D\big)^{-1}\big(B^\top \eta_\e+D^\top\z_\e+D^\top P_\e\si+\rho\big).\ea\right.\ee
Then
\bel{Mini-seq}u_\e(\cd)\deq \Th_\e(\cd)X_\e(\cd)+v_\e(\cd),\qq \e>0\ee
is a minimizing sequence of $u(\cd)\mapsto J(t,x;u(\cd))$:
\bel{Mini-seq-lim}\lim_{\e\to0}J(t,x;u_\e(\cd))=\inf_{u(\cd)\in\cU[t,T]}J(t,x;u(\cd))=V(t,x).\ee

\ms

\it Proof. \rm For any $\e>0$,  consider state equation (\ref{state}) and the
following cost functional:
\bel{} J_\e(t,x;u(\cd))\deq J(t,x;u(\cd))+\e\dbE\int_t^T|u(s)|^2ds.\ee
Denote the above problem by Problem $\hb{(SLQ)}_\e$ and the corresponding  value function by $V_\e(\cd\,,\cd)$.
By Corollary 4.8, $u_\e(\cd)$ defined by (\ref{Mini-seq}) is the unique optimal control of Problem $\hb{(SLQ)}_\e$
at $(t,x)\in[0,T)\times\dbR^n$. Note that
$$\ba{ll}
\ns\ds\e\dbE\int_t^T|u_\e(s)|^2ds=J_\e(t,x;u_\e(\cd))-J(t,x;u_\e(\cd))=V_\e(t,x)-J(t,x;u_\e(\cd))\\
\ns\ds\qq\qq\qq\q~\les V_\e(t,x)-V(t,x)\to 0\qq \hb{as}\q \e\to 0.\ea$$
Thus,
$$\lim_{\e\to0}J(t,x;u_\e(\cd))=\lim_{\e\to0}\bigg[V_\e(t,x)-\e\dbE\int_t^T|u_\e(s)|^2ds\bigg]=V(t,x).$$
The proof is completed. \endpf

\ms

Using the minimizing sequence constructed in Theorem 6.1, the open-loop solvability of Problem (SLQ)
can be characterized as follows.

\ms

\bf Theorem 6.2. \sl Let {\rm(H1)--(H2)} hold. Suppose
$u(\cd)\mapsto J^0(0,0;u(\cd))$ is convex. Let
$(t,x)\in[0,T)\times\dbR^n$ and $\{u_\e(\cd)\}_{\e>0}$ be the
sequence defined by {\rm(\ref{Mini-seq})}. Then the following
statements are equivalent:

\ms

{\rm(i)} Problem {\rm(SLQ)} is open-loop solvable at $(t,x)$;

\ms

{\rm(ii)} The sequence $\{u_\e(\cd)\}_{\e>0}$ admits a weakly
convergent subsequence;

\ms

{\rm(iii)} The sequence $\{u_\e(\cd)\}_{\e>0}$ admits a strongly
convergent subsequence.

\ms

In this case, the weak {\rm(}strong{\rm)} limit of any weakly
{\rm(}strongly{\rm)} convergent subsequence of
$\{u_\e(\cd)\}_{\e>0}$ is an open-loop optimal control of Problem {\rm(SLQ)}
at $(t,x)$.

\ms

\rm

To prove Theorem 6.2, we need the following lemma.

\ms

\bf Lemma 6.3. \sl Let $\cH$ be a Hilbert space with norm $\|\cd\|$ and $\th,
\th_n\in\cH$, $n=1,2,\cdots$.

\ms

{\rm(i)} If $\th_n \to \th$ weakly, then
$\ds\|\th\|\les\liminf_{n\to\i}\|\th_n\|$.

\ms

{\rm(ii)} $\th_n \to \th$ strongly if and only if
$$\|\th_n\|\to\|\th\|\qq \hb{and} \qq\th_n\to\th \hb{\q weakly}.$$

\ms

\it Proof. \rm (i) By the Hahn-Banach theorem, we can choose a $\th^*\in\cH$ with $\|\th^*\|=1$ such that
$\langle \th^*,\th\rangle=\|\th\|$. Thus,
$$\|\th\|=\langle\th^*,\th\rangle=\lim_{n\to\i}\langle \th^*,\th_n\rangle\les\liminf_{n\to\i}\|\th_n\|.$$

\ms

(ii) The necessity is obvious. Now suppose $\|\th_n\|\to\|\th\|$ and $\th_n\to\th$ weakly. Then
$$\|\th_n-\th\|=\|\th_n\|^2-2\lan\th,\th_n\ran+\|\th\|^2\to 0\qq\hb{as}\q n\to\i.$$
This completes the proof.
\endpf

\ms

\it Proof of Theorem {\rm6.2.} \rm (i) $\Ra$ (ii), and (i) $\Ra$ (iii).
Let $v^*(\cd)$ be an open-loop optimal control of Problem (SLQ) at $(t,x)$. By
Corollary 4.8, for any $\e>0$, $u_\e(\cd)$ defined by
(\ref{Mini-seq}) is the unique optimal control of Problem
${\rm(SLQ)}_\e$ at $(t,x)$ and
\bel{w-1}V_\e(t,x)=J_\e(t,x;u_\e(\cd))\ges V(t,x)+\e\dbE\int_t^T|u_\e(s)|^2ds.\ee
Also, we have
\bel{w-2}V_\e(t,x)\les J_\e(t,x;v^*(\cd))=V(t,x)+\e\dbE\int_t^T|v^*(s)|^2ds.\ee
Combining (\ref{w-1})--(\ref{w-2}), we have
\bel{w-3}\dbE\int_t^T|u_\e(s)|^2ds\les{V_\e(t,x)-V(t,x)\over\e}\les\dbE\int_t^T|v^*(s)|^2ds,\qq\forall \e>0.\ee
Thus, $\{u_\e(\cd)\}_{\e>0}$ is bounded in the Hilbert space $\cU[t,T]=L_\dbF^2(t,T;\dbR^m)$ and hence admits
a weakly convergent subsequence $\{u_{\e_k}(\cd)\}_{k\ges1}$.
Let $u^*(\cd)$ be the weak limit of $\{u_{\e_k}(\cd)\}_{k\ges1}$. Since $u(\cd)\mapsto J(t,x;u(\cd))$ is convex and continuous,
it is hence sequentially weakly lower semi-continuous. Thus (noting (\ref{Mini-seq-lim})),
$$V(t,x)\les J(t,x;u^*(\cd))\les \liminf_{k\to\i}J(t,x;u_{\e_k}(\cd))=V(t,x),$$
which implies that $u^*(\cd)$ is also an open-loop optimal control of Problem (SLQ) at $(t,x)$.
Now replacing $v^*(\cd)$ with $u^*(\cd)$ in (\ref{w-3}), we have
\bel{w-4}\dbE\int_t^T|u_\e(s)|^2ds\les\dbE\int_t^T|u^*(s)|^2ds,\qq\forall \e>0.\ee
Also, by Lemma 6.3, part (i),
\bel{w-5}\dbE\int_t^T|u^*(s)|^2ds\les\liminf_{k\to\i}\dbE\int_t^T|u_{\e_k}(s)|^2ds.\ee
Combining (\ref{w-4})--(\ref{w-5}), we see that
$$\dbE\int_t^T|u^*(s)|^2ds=\lim_{k\to\i}\dbE\int_t^T|u_{\e_k}(s)|^2ds.$$
Then it follows from Lemma 6.3, part (ii), that
$\{u_{\e_k}(\cd)\}_{k\ges1}$ converges to $u^*(\cd)$ strongly.
\ms

(iii) $\Ra$ (ii) is obvious.

\ms

(ii) $\Ra$ (i). Let $\{u_{\e_k}(\cd)\}_{k\ges1}$ be a weakly convergent subsequence of $\{u_\e(\cd)\}_{\e>0}$ with weak
limit $u^*(\cd)$. Then $\{u_{\e_k}(\cd)\}_{k\ges1}$ is bounded in $\cU[t,T]=L_\dbF^2(t,T;\dbR^m)$. For any $u(\cd)\in\cU[t,T]$,
we have
\bel{w-6}J(t,x;u(\cd))+\e_k\dbE\int_t^T|u(s)|^2ds\ges V_{\e_k}(t,x)=J(t,x;u_{\e_k}(\cd))+\e_k\dbE\int_t^T|u_{\e_k}(s)|^2ds.\ee
Note that $u(\cd)\mapsto J(t,x;u(\cd))$ is sequentially weakly lower semi-continuous. Letting $k\to\i$ in (\ref{w-6}), we obtain
$$J(t,x;u^*(\cd))\les \liminf_{k\to\i}J(t,x;u_{\e_k}(\cd))\les J(t,x;u(\cd)),\qq\forall u(\cd)\in\cU[t,T].$$
Hence, $u^*(\cd)$  is an open-loop optimal control of Problem (SLQ) at $(t,x)$.
\endpf

\ms

From the proof of Theorem 6.2, we see that the open-loop solvability
of Problem {\rm(SLQ)} at $(t,x)$ is also equivalent the
$L^2$-boundedness of $\{u_\e(\cd)\}_{\e>0}$.
In particular, the open-loop solvability of Problem $\hb{(SLQ)}^0$ at $(t,x)$ is equivalent the
$L^2$-boundedness of $\{\Th_\e(\cd)X_\e(\cd)\}_{\e>0}$ with $X_\e(\cd)$ being the solution of
\bel{}\left\{\2n\ba{ll}
\ns\ds dX_\e(s)=(A+B\Th_\e)X_\e ds+(C+D\Th_\e)X_\e dW(s),\qq s\in[t,T], \\
\ns\ds X_\e(t)=x.\ea\right.\ee
Since the
$L_\dbF^2(\O;C([t,T];\dbR^n))$-norm of $X_\e(\cd)$ is dominated by the
$L^2$-norm of $\Th_\e(\cd)$, we conjecture that the
$L^2$-boundedness of $\{\Th_\e(\cd)\}_{\e>0}$ will lead to the
open-loop solvability of Problem $\hb{(SLQ)}^0$ at $(t,x)$. Actually, we have the
following result.

\ms

\bf Proposition 6.4. \sl Let {\rm(H1)--(H2)} hold. Suppose
$u(\cd)\mapsto J^0(0,0;u(\cd))$ is convex, and let
$\{\Th_\e(\cd)\}_{\e>0}$ be the sequence defined by
{\rm(\ref{Th-e})}. If
\bel{Mini-prop-0}\sup_{\e>0}\int_0^T|\Th_\e(s)|^2ds<\i,\ee
then the Riccati equation {\rm(\ref{Riccati})} admits a
regular solution $P(\cd)\in C([0,T];\dbS^n)$. Consequently, Problem
${\rm(SLQ)}^0$ is closed-loop solvable.

\ms

\it Proof. \rm For any $x\in\dbR^n$ and $\e>0$, let $X_\e(\cd)$ be the solution of
\bel{}\left\{\2n\ba{ll}
\ns\ds dX_\e(s)=\big[A(s)+B(s)\Th_\e(s)\big]X_\e(s)ds+\big[C(s)+D(s)\Th_\e(s)\big]X_\e(s)dW(s),\qq s\in[0,T], \\
\ns\ds X_\e(0)=x.\ea\right.\ee
By It\^o's formula, we have
$$\ba{ll}
\ns\ds\dbE|X_\e(t)|^2=|x|^2+\dbE\int_0^t\Big\{\big|[C(s)+D(s)\Th_\e(s)]X_\e(s)\big|^2+2\lan[A(s)+B(s)\Th_\e(s)]X_\e(s),X_\e(s)\ran\Big\}ds\\
\ns\ds\qq\qq~\2n\les|x|^2+\int_0^t\[\big|C(s)+D(s)\Th_\e(s)\big|^2+2\big|A(s)+B(s)\Th_\e(s)\big|\]\dbE|X_\e(s)|^2ds,\qq\forall t\in[0,T].
\ea$$
Thus, by Gronwall's inequality,
\bel{}\ba{ll}
\ns\ds\dbE|X_\e(t)|^2\les|x|^2\exp\left\{\int_0^T\[|C(s)+D(s)\Th_\e(s)|^2+2|A(s)+B(s)\Th_\e(s)|\]ds\right\}\\
\ns\ds\qq\qq~\2n\les|x|^2\exp\left\{K\left(1+\int_0^T|\Th_\e(s)|^2ds\right)\right\},\qq\forall t\in[0,T],\ea\ee
where $K>0$ is some constant depending only on $A(\cd),B(\cd),C(\cd),D(\cd)$. Hence,
\bel{}\ba{ll}
\ns\ds\dbE\int_0^T|\Th_\e(s)X_\e(s)|^2ds\les\int_0^T|\Th_\e(s)|^2\dbE|X_\e(s)|^2ds\\
\ns\ds\qq\qq\qq\qq\qq\1n\les|x|^2\exp\left\{K\left(1+\int_0^T|\Th_\e(s)|^2ds\right)\right\}\int_0^T|\Th_\e(s)|^2ds,\ea\ee
which, together with (\ref{Mini-prop-0}), implies the
$L^2$-boundedness of $\{\Th_\e(s)X_\e(s)\}_{\e>0}$. Thus, by Theorem
6.2, Problem $\hb{(SLQ)}^0$ is open-loop solvable at $t=0$, and by Theorem 5.3,
Problem $\hb{(SLQ)}^0$ is finite. Now let
$P:[0,T]\to\dbS^n$ such that {\rm(\ref{V0})} holds. Then by Theorem 5.3,
\bel{prop7.4-1}R+D^\top PD\ges0,\qq\ae\ee
Let $\{\Th_{\e_k}(\cd)\}$ be a weakly convergent subsequence of
$\{\Th_\e(\cd)\}$ with weak limit $\Th(\cd)$. Since
$$R+\e I+D^\top P_\e D\to R+D^\top PD\qq\hb{as}\q \e\to0$$
and $\{R(\cd)+\e I+D(\cd)^\top P_\e(\cd)D(\cd)\}_{0<\e\les1}$ is uniformly bounded,
we have
$$B^\top P_{\e_k}+D^\top P_{\e_k} C+S=-(R+\e_k I+D^\top P_{\e_k} D)\Th_{\e_k}\to-(R+D^\top PD)\Th \qq\hb{weakly in } L^2.$$
Also, note that
$$B^\top P_{\e_k}+D^\top P_{\e_k} C+S\to B^\top P+D^\top P C+S \qq\hb{strongly in } L^2.$$
Thus,
$$-(R+D^\top PD)\Th=B^\top P+D^\top P C+S.$$
This implies
\bel{prop7.4-2}\cR\big(B^\top P+D^\top PC+S\big)\subseteq\cR\big(R+D^\top PD\big),\qq\ae\ee
Since
$$(R+D^\top PD)^\dag(B^\top P+D^\top PC+S)=-(R+D^\top PD\big)^\dag(R+D^\top PD\big)\Th,$$
and $(R+D^\top PD\big)^\dag(R+D^\top PD\big)$ is an orthogonal projection, we have
\bel{prop7.4-3}(R+D^\top PD)^\dag(B^\top P+D^\top PC+S)\in L^2(0,T;\dbR^{m\times n}),\ee
and
$$\Th=-(R+D^\top PD)^\dag(B^\top P+D^\top PC+S)+\big[I-(R+D^\top PD)^\dag(R+D^\top PD)\big]\Pi$$
for some $\Pi(\cd)\in L^2(0,T;\dbR^{m\times n})$. Finally, letting $k\to\i$, we have
$$\ba{ll}
\ns\ds P(t)=\lim_{k\to\i}P_{\e_k}(t)=G+\lim_{k\to\i}\int_t^T\[P_{\e_k} A+A^\top P_{\e_k}+C^\top P_{\e_k} C+Q
+(P_{\e_k} B+C^\top P_{\e_k} D+S^\top)\Th_{\e_k}\]ds\\
\ns\ds\qq\,=G+\int_t^T\[PA+A^\top P+C^\top PC+Q+(PB+C^\top PD+S^\top)\Th\]ds\\
\ns\ds\qq\,=G+\int_t^T\[PA\1n+\1nA^\top P\1n+\1nC^\top PC\1n+\1nQ\1n
-\1n(PB\1n+\1nC^\top PD\1n+\1nS^\top)(R\1n+\1nD^\top PD)^\dag(B^\top P\1n+\1nD^\top PC\1n+\1nS)\]ds,
\ea$$
which, together with (\ref{prop7.4-1})--(\ref{prop7.4-3}), implies that $P(\cd)$ is a regular solution of (\ref{Riccati}).
\endpf

\section{An example}

In this section we re-exam Example 2.2 to illustrate some results we
obtained. In this example, the stochastic LQ problem admits a {\it
continuous} open-loop optimal control at all
$(t,x)\in[0,T)\times\dbR^n$, hence it is open-loop solvable, while
the value function is {\it not} continuous in $t$; the corresponding
Riccati equation has a unique solution $P(\cd)$, which does {\it
not} satisfy the range condition (\ref{regular-1}) and therefore is
not regular. Therefore, the problem is {\it not} closed-loop
solvable on any $[t,T]$. This example also tells us that the
necessity part in Theorem 4.2 does not hold.

\ms

Recall the following Problem $\hb{(SLQ)}^0$ with one-dimensional state equation:
\bel{ex1-1}\left\{\2n\ba{ll}
\ns\ds dX(s)=\big[u_1(s)+u_2(s)\big]ds+\big[u_1(s)-u_2(s)\big]dW(s),\qq s\in[t,1], \\
\ns\ds X(t)=x,\ea\right.\ee
and cost functional
\bel{ex1-2}J^0(t,x;u(\cd))=\dbE X(1)^2.\ee
In this example, $u(\cd)=(u_1(\cd),u_2(\cd))^\top$ and
$$\left\{\2n\ba{ll}
\ns\ds A=0,\q B=(1,1),\q C=0,\q D=(1,-1),\\
\ns\ds G=1,\q Q=0,\q S=(0,0)^\top,\q R=\begin{pmatrix}0&0\\0&0\end{pmatrix}.\ea\right.$$
The corresponding Riccati equation reads
\bel{ex1-3}\left\{\2n\ba{ll}
\ns\ds\dot P=P^2(1,1)\begin{pmatrix}P&-P\\-P&P\end{pmatrix}^\dag\begin{pmatrix}1\\1\end{pmatrix}
={P\over4}(1,1)\begin{pmatrix}1&-1\\-1&1\end{pmatrix}\begin{pmatrix}1\\1\end{pmatrix}=0,\\
\ns\ds P(1)=1.\ea\right.\ee
Obviously, (\ref{ex1-3}) has a unique solution $P(\cd)\equiv1$, and
$$\ba{ll}
\ns\ds \cR\big(B(s)^\top P(s)+D(s)^\top P(s)C(s)+S(s)\big)=\cR\big((1,1)^\top\big)=\left\{(a,a)^\top: a\in\dbR\right\},\\
\ns\ds \cR\big(R(s)+D(s)^\top P(s)D(s)\big)=\cR\left(\begin{pmatrix}1&-1\\-1&1\end{pmatrix}\right)=\left\{(a,-a)^\top: a\in\dbR\right\}.\ea$$
Thus, the range condition (\ref{regular-1}) does not hold and hence $P(\cd)$ is not regular.
By Theorem 4.4, the problem is not closed-loop solvable on any $[t,1]$.

\ms

Now for any $\e>0$, consider state equation (\ref{ex1-1}) and the cost functional
\bel{ex1-Je}J_\e^0(t,x;u(\cd))=\dbE \left[X(1)^2+\e\int_t^T|u(s)|^2ds\right].\ee
The Riccati equation for the above problem reads
\bel{ex1-4}\left\{\2n\ba{ll}
\ns\ds\dot P_\e=P_\e^2(1,1)\begin{pmatrix}\e+P_\e&-P_\e\\-P_\e&\e+P_\e\end{pmatrix}^{-1}\begin{pmatrix}1\\1\end{pmatrix}
={2\over\e}P_\e^2,\\
\ns\ds P_\e(1)=1,\ea\right.\ee
whose solution is given by
\bel{ex1-5}P_\e(t)={\e\over\e+2-2t},\qq t\in[0,1].\ee
Letting $\e\to 0$, we have
\bel{ex1-6}P_0(t)\deq\lim_{\e\to0}P_\e(t)=\left\{\2n\ba{ll}
\ns\ds 0,\qq 0\les t<1,\\
\ns\ds 1,\qq t=1.\ea\right.\ee
Thus, by Theorem 5.3, the original Problem $\hb{(SLQ)}^0$ is finite with value function
\bel{ex1-7}V^0(t,x)=0,\q 0\les t<1;\qq V^0(1,x)=x^2,\qq\forall x\in\dbR.\ee
Next, set
$$\Th_\e\deq-(R+\e I+D^\top P_\e D)^{-1}(B^\top P_\e+D^\top P_\e C+S)=-{P_\e\over\e}\begin{pmatrix}1\\1\end{pmatrix}.$$
Then the solution of
\bel{ex1-8}\left\{\2n\ba{ll}
\ns\ds dX_\e(s)=\big[A(s)+B(s)\Th_\e(s)\big]X_\e(s)ds+\big[C(s)+D(s)\Th_\e(s)\big]X_\e(s)dW(s)\\
\ns\ds\qq\q~\2n=-{2P_\e\over\e}X_\e(s)ds,\qq s\in[t,T], \\
\ns\ds X_\e(t)=x\ea\right.\ee
is given by
\bel{ex1-9}X_\e(s)=x\exp\left\{-\int_t^s{2P_\e(r)\over\e}dr\right\}={\e+2-2s\over\e+2-2t}x,\qq t\les s\les1,\ee
and hence
\bel{ex1-10}u_\e(s)\deq\Th_\e(s)X_\e(s)=-\left({x\over\e+2-2t},{x\over\e+2-2t}\right)^\top,\qq t\les s\les1.\ee
Note that for $t\in[0,1)$,
$$u_\e(\cd)\to -\left({x\over2-2t},{x\over2-2t}\right)^\top\q\hb{in } L^2\hb{ as }\e\to0.$$
Thus, by Theorem 6.2, the original Problem $\hb{(SLQ)}^0$ is open-loop solvable at any
$(t,x)\in[0,T)\times\dbR$ with an open-loop optimal control
\bel{ex1-11}u^*_{(t,x)}(s)=-\left({x\over2-2t},{x\over2-2t}\right)^\top,\qq t\les s\les1,\ee
which is continuous in $s\in[t,1]$.


\end{document}